\documentclass[conference]{IEEEtran}
\pdfoutput=1

\usepackage{dsfont}
\usepackage[T1]{fontenc}
\usepackage{graphicx}
\usepackage{soul}
\usepackage{color}
\usepackage[utf8x]{inputenc}
\usepackage{amsfonts}
\usepackage{amsmath}
\usepackage{setspace}
\usepackage{amssymb}
\usepackage{booktabs,graphicx}
\usepackage[noadjust]{cite}
\usepackage{bbm}
\usepackage{wasysym}
\usepackage{textcomp}
\IEEEoverridecommandlockouts

\DeclareMathAlphabet{\mymathbb}{U}{bbold}{m}{n}

\DeclareMathAlphabet{\mymathbb}{U}{BOONDOX-ds}{m}{n}

\newtheorem{defn}{Definition}
\newtheorem{rmk}{Remark}

\def\L2{{\cal L}_2}
\def\L2e{{\cal L}_{2e}}

\def\lef[{\left[\begin{array}}
	\def\rig]{\end{array}\right]}
\def\begrem{\begin{remark}}
	\def\endrem{\end{remark}}
\def\begequarr{\begin{eqnarray}}
\def\endequarr{\end{eqnarray}}
\def\begequarrs{\begin{eqnarray*}}
	\def\endequarrs{\end{eqnarray*}}
\def\begarr{\begin{array}}
	\def\endarr{\end{array}}
\def\begequ{\begin{equation}}
\def\endequ{\end{equation}}

\def\begdes{\begin{description}}
	\def\enddes{\end{description}}
\def\begenu{\begin{enumerate}}
	\def\begite{\begin{itemize}}
		\def\endite{\end{itemize}}
	\def\endenu{\end{enumerate}}

\usepackage{color}
\newcommand{\yo}[1]{{\color{black} #1}}

\newcommand{\jdw}[1]{{\color{black} #1}}
\newcommand{\icl}[1]{{\color{black} #1}}

\newcommand{\rmspace}{\vspace{-6mm}}

\makeatletter
\IEEEtriggercmd{\reset@font\normalfont\footnotesize}
\makeatother
\IEEEtriggeratref{1}

\begin{document}

\title{ A Review of Reduced-Order Models for Microgrids: Simplifications vs Accuracy}


\author{	
	\IEEEauthorblockN{Yemi Ojo, Jeremy Watson, and Ioannis Lestas}
	\IEEEauthorblockA{Department of Engineering, University of Cambridge \\ Cambridge, United Kingdom \\ yo259@cam.ac.uk, jdw69@cam.ac.uk, icl20@cam.ac.uk}
}

	\maketitle

	\begin{abstract}
	
Inverter-based  microgrids are an important technology for sustainable electrical power systems and typically use droop-controlled grid-forming inverters to interface distributed energy resources to the network and control the voltage and frequency.
	\jdw{\icl{Ensuring stability of such microgrids is a key issue, which requires the use of appropriate models for  analysis and control system design.}}
\yo{Full-order detailed models \icl{can be more difficult to analyze  and increase computational complexity,}
hence a number \icl{of} reduced-order models have been proposed in the literature which present various trade-offs between accuracy and complexity.}
However, the\jdw{se} simplifications present the risk of \jdw{failing to} adequately capture important dynamics of the microgrid.
\jdw{Therefore, there is a need for} a comprehensive review and assessment of their relative quality, 
\icl{which is something that
has not been systematically carried out thus far in the literature
and we aim to address in this paper.}
In particular, we review various inverter-based microgrid reduced-order models and investigate
\icl{the accuracy of their predictions for stability via
a comparison with
a corresponding detailed average model.}
Our study shows that \icl{the simplifications reduced order models rely upon can 
affect their accuracy 
in various regimes of} 
the line R/X ratios,  and that inappropriate model choices can result in substantially inaccurate stability results. Finally, we present recommendations \jdw{on} the use \icl{of reduced order} models for the stability \icl{analysis} of microgrids.

	\end{abstract}

\begin{IEEEkeywords}
	Microgrids, droop control, reduced-order model, small-signal stability.
\end{IEEEkeywords}

	\section{INTRODUCTION}%
	\label{intro}

Increasing environmental \icl{concerns}, advancement in renewable energy technologies and the continuous rise in the global energy demand bring about the exploration of renewable energy sources. Most of these sources are small-scale distributed energy generation (DG) units and will form microgrids, which are the interconnection of DG units, loads, and energy storage systems into a controllable system. Microgrids can be operated in grid-connected or autonomous mode, and autonomous microgrids are an especially effective solution for remote areas where the main power grid is not accessible \cite{hatziargyriou2007, lasseter2002,katiraei2008}.
Most renewable energy sources cannot be effectively utilised without power conditioners \cite{hill2012battery,singh2013solar, blaabjerg2011power, bresesti2007hvdc, barakati2005}, and power converters are  becoming increasingly common as interfaces that provide the power conditioning capability. They now appear in higher power ratings and are therefore an integral part of microgrid dynamics.

There are new challenges associated with microgrids. With the increasing proportion of DG units, especially in the autonomous mode of operation, there is a need for individual inverters to achieve power sharing. Power sharing is the ability of microgrid inverters to optimally allocate their power outputs to meet the load demand in the network, while achieving a desired steady-state (characterised by the steady-state frequency and voltage). It is desirable to achieve this via the local measurements (voltages and currents) provided to the inverters' local controllers. The droop control policy offers this feature and is widely used \cite{chandorkar1993, pogaku2007, vorobev2018, Guerrero2014, PogakuBarklund2008, yemi2017, yemimo2019,yemijeremy2019, sao2005,ysunGuerrero2017}. Local measurements are used to adjust the inverters' frequency and voltage setpoints based on the corresponding demand of active and reactive power respectively. The droop gains are chosen to specify the active and reactive power sharing ratio among inverters, however this choice strongly affects the dynamics and  stability of microgrids \cite{pogaku2007, vorobev2018, Guerrero2014, Guerrero2015, coelho2002, coelho2000}.
The nonlinearity of this control scheme increases the complexity of the closed loop microgrid dynamics, and severe stability issues have been noted in previous studies \cite{greengu2017reduced, boroyevich2013intergrid, flourentzou2009vsc}. There is hence a need for stability issues to be accurately assessed via appropriate mathematical models.
	
Accurate mathematical representation of inverters is fundamental for the modelling and stability assessment of inverter-based microgrids. Inverters are usually controlled \icl{via pulse-width modulation (PWM) schemes} \cite{blaabjergXiaolong2019, Rocabert2012, nabae1981, holmes2003}. The \icl{discontinuity of trajectories in this scheme
\yo{increases the complexity of inverter dynamics}
 and a common approach to address this is the use of the averaging theory,} which allows inverters to be represented as continuous time systems \cite{ mollerstedt2000out, blaabjergharneforspassivity2015, turner2012case, pogaku2007, jeremyyemi2019}. 
  Such models are used \icl{either} \jdw{for} frequency domain \jdw{analysis} \cite{turner2012case, blaabjergwangpassivity2017, blaabjergwang2017unified, harnefors2017, blaabjergwang2014modeling, paice2000rail}, or \jdw{analysis in the} time domain \cite{pogaku2007, greengu2017reduced, greenkroutikova2007state}. The latter allows for both small-signal and large signal models to be formulated,
  \icl{and provides} additional insight into dynamic interactions \cite{pogaku2007}.

Models used to describe the behaviour of inverter-based microgrids consist of the inverter and power lines dynamics, and typically employ the averaging theory-based inverter models due to their versatility and easy application.
 \icl{A detailed} average model which consists of all the internal states of an inverter as well as the network dynamics was developed in \cite{pogaku2007}.
It has been extensively \icl{used for} stability assessment of microgrids of different sizes and configurations, and \icl{is known} \icl{to give} 
reliable stability  results.

\icl{In order to facilitate analytical results at the network level various reduced-order models have been used in the literature.}
\jdw{These} generally assume that inverters are controllable voltage sources (i.e. the internal controllers of inverters are fast and track their reference signals). The line dynamics are either considered in full or approximated to various degrees.
\icl{Such simplified}
reduced models can be broadly classified into three models: the 5th-order model, \icl{the} conventional 3rd-order model, and \icl{a high-fidelity} 3rd-order model.
The 5th-order model takes into account the dynamical behaviour of the lines and represents the inverter as a first-order voltage source \cite{Guerrero2014, Guerrero2015, vorobev2018}.
The conventional 3rd-order model uses a traditional power system assumption of a distinct timescale separation which allows the line dynamics to be neglected and the lines modelled statically \cite{kundur94, machowski08, padiyar1996, anderson02, glover11}. In addition, lossless lines are assumed and the inverter is modelled as a first-order voltage source \cite{schiffer2014, guerreroMatas2004, sao2005, guerreroMatas2006, coelho2002, coelho2000, chandorkar1993}. The Kuramoto model in \cite{simpson2013} additionally assumes a fixed voltage magnitude.
However, due to the low inertia of inverters, a timescale separation approach is less valid for inverter-based microgrids. Fast dynamics (e.g. the fast timescale associated with the small $X/R$ ratio of power lines) can significantly influence the dynamics of slower modes \cite{federicoDofler2018, Guerrero2015, vorobev2018}.
The high-fidelity 3rd-order model  \cite{vorobev2018,vorobev2017} sought to address this issue. It incorporates a first-order Taylor \icl{series} approximation of the line dynamics, and again considers the inverter as a first-order voltage source.

\icl{The simplicity of reduced order models can be of value analytically and computationally.  Nevertheless, they raise the risk of failing to capture important aspects of the dynamics of microgrids \cite{hatziargyriou2014microgrids}.
The sensitivity of microgrid stability to line dynamics, and the strong coupling between voltage and frequency dynamics 
pose the need for a thorough investigation of the accuracy of such models. In particular, the performance of these simplified models and the regimes where they are suitable for stability assessment must be clarified. 

This is a main aim of this paper and its contributions can be summarised as follows.}

	\begin{itemize}
		\item We investigate the accuracy
		\yo{of reduced order \icl{average} models for  microgrids  used in the literature} \icl{which include}
		the conventional 3rd-order, \icl{the} high-fidelity 3rd-order, and \icl{the} electromagnetic 5th-order models. This is conducted via a comparison of \icl{the regime they predict stability to that of a detailed} average model.
		\item
\icl{We show that various simplifications upon which the reduced models rely  have a significant impact on their accuracy in different regimes of the line $R/X$ ratios. Hence inappropriate model choices can result in substantially inaccurate stability results.
		\item	Based on our findings we provide recommendations for future studies on
appropriate 
models to be used for microgrid stability assessment for different values of the line $R/X$ ratios.}	
	\end{itemize}

	The paper is structured as follows. Section \ref{nota} gives the notation and definitions used in the paper. An overview  of the work  is presented in section \ref{over}. Section \ref{models} presents the detailed average model and the three reduced-order nonlinear models to be considered. Their small-signal models are presented in section \ref{smallsignal}, and we evaluate the correctness of stability regions  of the reduced-order models  in section \ref{accuracy}.	Section \ref{disc} discusses the implications of our study, and \icl{conclusions of the paper are provided} in section \ref{concl}.

	\section{Notation and Definitions}
	\label{nota}
	Let $\mathbb{R}_{\geq0}:=\{x\in\mathbb{R}|x\geq0\}$,  $\mathbb{R}_{>0}:=\{x\in\mathbb{R}|x>0\}$ and $\mathbb{S}\in[0, 2\pi)$.
	Given $n, p \in \mathbb{N^*}$,
$\mathds{1}_n$ denotes the column vector of ones with length $n$,
$\mathbf{e}=\scriptsize{[1~~0]^{\top}}$,
	 $\mathbf{I}_{n}$ as the $n \times n$ identity matrix, $\mymathbb{0}_{n\times p}$ the $n \times p$ matrix of all zeros, $J=\scriptsize{\begin{bmatrix}
	0&1\\-1&0	\end{bmatrix}}$,
and diag($a_i$), $i=1\ldots n$, an $n\times n$ diagonal matrix with diagonal entries $a_i$.

\begin{defn}[Symmetric AC three-phase signals] 
	\label{def1}
	$x_{abc}:\mathbb{R}_{>0} \to \mathbb{R}^3$ denotes a symmetric three-phase AC signal of the form
	 $x_{abc}(t)=\text{col}(\hat{X}(t)\sin(\theta(t)),\hat{X}(t)\sin(\theta(t)-\frac{2\pi}{3}),\hat{X}(t)\sin(\theta(t)+\frac{2\pi}{3}))$
	where $\theta(t)\in\mathbb{R}$ is the angle that satisfies
	\begin{equation}
	\label{theta}
	\dot \theta(t)=\omega(t)
	\end{equation}
	and $\hat{X}(t)$, $\omega(t)$  are respectively the amplitude and frequency which are in general time varying variables that take values in $\mathbb{R}_{>0}$. Note also that  $\mathds{1}_3^{\top}x_{abc}=0$.
\end{defn}
The time argument $t$ will often be omitted in the presentation for convenience in the notation.

\begin{defn}
	\label{def4}
	The Park transformation matrix for a given $\theta$
	is given by
	\begin{equation}\scriptsize
	\label{park1}
	T(\theta)=\sqrt{\frac{2}{3}}\begin{bmatrix}  \sin(\theta) &\sin(\theta-\frac{2\pi}{3}) &\sin(\theta+\frac{2\pi}{3}) \\
	\cos(\theta) &\cos(\theta-\frac{2\pi}{3}) &\cos(\theta+\frac{2\pi}{3}) \\ \frac{1}{\sqrt{2}} &\frac{1}{\sqrt{2}} &\frac{1}{\sqrt{2}} \end{bmatrix}
	\end{equation}
\end{defn}

	\begin{defn}[Local Direct-Quadrature Coordinates] \label{def3a}
		\label{def3}
		The representation of a signal $x_{abc}$ (Definition \ref{def1}) in its Local Direct-Quadrature-Zero coordinates (denoted by $x_{dq0}$)
\icl{is given by}
		\begin{equation}
		\label{eqn5a}
		x_{dq0}(t)=T(\theta(t))x_{abc}(t).
		\end{equation}		
		where $T(\theta(t))$ is as in \eqref{park1} with $\theta(t)$ the angle of signal $x_{abc}(t)$ as in Definition \ref{def1}.
		Note that $x_{abc}$ can be obtained from a given  $x_{dq0}$ using \eqref{eqn5c}.
		\begin{equation}
		\label{eqn5c}
		x_{abc}(t)=T^{-1}(\theta(t))x_{dq0}(t)
		\end{equation}
	\end{defn}

	\begin{defn}[Synchronous Direct-Quadrature Coordinates] \label{DQ0}
		The representation of a signal $x_{abc}$ (Definition \ref{def1}) in its Synchronous Direct-Quadrature-Zero coordinates (denoted by $x_{DQ0}$)
\icl{is given by}
		\begin{equation}
		\label{eqn5b}
		x_{DQ0}(t)=T(\theta_0)x_{abc}(t).
		\end{equation}	
		where $T(\theta_0)$ \icl{is as in \eqref{park1} and the angle $\theta_0$ 
satisfies}
		\begin{equation}
		\label{theta0}
		\dot \theta_0=\omega_0
		\end{equation}
		where $\omega_0\in\mathbb{R}_{>0}$ is \icl{a constant} synchronous 
frequency,
		and the signal $x_{abc}(t)$ \yo{\icl{is as in Definition \ref{def1}}}.
		Note also that $x_{abc}$ can be obtained from a given  $x_{DQ0}$ using \icl{the relation} 
		\begin{equation}
		\label{eqn5d}
		x_{abc}(t)=T^{-1}(\theta_0)x_{DQ0}(t)
		\end{equation}
	\end{defn}

\begin{rmk}\label{rem:dq}
	Since $x_{abc}$ is symmetric by Definition \ref{def1}, the third component of $x_{dq0}$ and $x_{DQ0}$ is zero. Hence $x_{dq}, x_{DQ}$ 
	refer to the first two entries of $x_{dq0}$ and $x_{DQ0}$ respectively.
\end{rmk}

\begin{rmk}[dq to DQ transformation]
	\label{def8}
	As a consequence of Definition \ref{def3a} and \ref{DQ0}, the relation between the quantities in the local $dq$ and synchronous $DQ$ frames is given by
	\begin{equation}
	\label{alg4}
	\begin{split}
	x_{DQ}(t)=&\mathcal{T}(\delta(t))x_{dq}(t),\\
	\mathcal{T}(\delta(t)) =&\scriptsize{\begin{bmatrix}
		\cos \delta(t) &-\sin \delta(t)\\ \sin \delta(t) & \cos\delta(t)
		\end{bmatrix}},\\
	\dot\delta(t)=&\omega(t)-\omega_{0}.
	\end{split}
	\end{equation}
\icl{Note that $\delta(t)$ is the angle between the two rotating frames of reference used in Definitions \ref{def3}, \ref{DQ0}, respectively. Also $\mathcal{T}(\delta(t))$ is a rotation matrix that satisfies} the properties: \mbox{$\mathcal{T}^{-1}(\delta)=\mathcal{T}^{\top}(\delta)$}, \mbox{$\frac{\partial \mathcal{T}(\delta)}{\partial\delta}=-J\mathcal{T}(\delta)$}, and
	${J}\mathcal{T}(\delta)=\mathcal{T}(\delta){J}$.
\end{rmk}		

	\section{\icl{Overview}}
	\label{over}
\icl{The aim of this study is to provide a rigorous review and a comparison of the stability properties of various reduced order models 
for interconnections of inverters with a grid forming role, i.e. inverters that aim to form an autonomous microgrid where they control its frequency and voltage.}



\begin{figure}[ht]
	\begin{center}
		\includegraphics[width=0.5\textwidth]{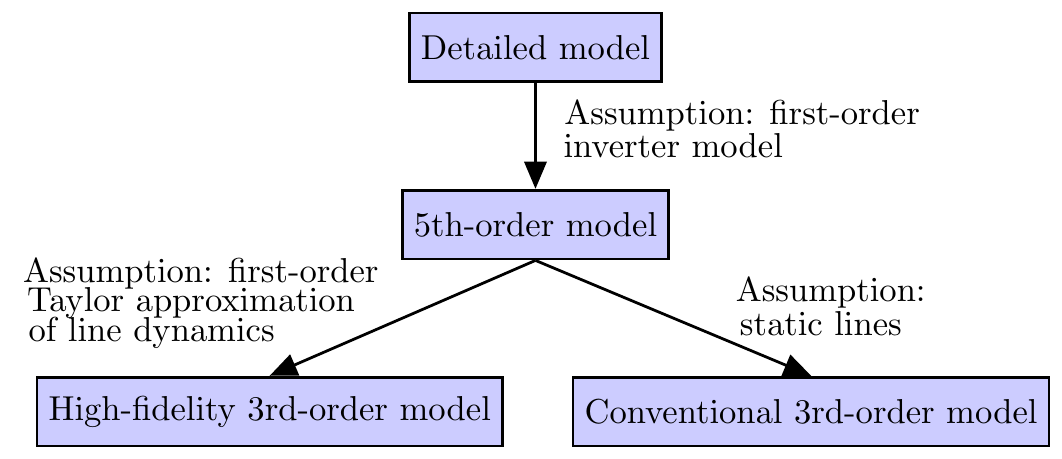}
		\vspace{-2mm}
		\caption{Detailed average model and  reduced-order models with associated assumptions.}
			\vspace{-4mm}
		\label{figflow}
	\end{center}
\end{figure}

\jdw{
\icl{Our study focuses on average models and as a benchmark for evaluating the accuracy of the reduced order models we will consider the detailed average model developed in \cite{pogaku2007} which has been extensively verified in the literature. The reduced order models we will consider are outlined below:}}
\begin{itemize}
\icl{\item {\em 5th-order model}: this takes into account the dynamical behaviour of the lines and represents the inverter as a first-order voltage source \cite{Guerrero2014, Guerrero2015, vorobev2018}.}
	\item \icl{{\em{Conventional 3rd-order model}}: this} uses a traditional power system assumption of a \icl{timescale separation between line and bus dynamics. This} allows the line dynamics to be neglected and the lines \icl{to be} modelled statically \cite{kundur94, machowski08, padiyar1996, anderson02, glover11}. In addition, the inverter is modelled as a first-order voltage source \cite{schiffer2014, vorobev2018}.
	\item \icl{{\em High-fidelity 3rd-order model:} this again} considers the inverter as a first-order voltage source, and    improves on the conventional 3rd-order model by  modelling the lines with a first-order Taylor approximation \cite{vorobev2018,vorobev2017} \icl{of their dynamics}.
\end{itemize}
Fig. \ref{figflow} \icl{summarises the connection between the detailed average model and the reduced-order models, and the assumptions the latter rely upon. A full description of these models will be provided in} section \ref{models}.

\subsection{\icl{System description}}
\begin{figure}[h]
	  	\begin{center}
	  		\includegraphics[width=0.3\textwidth]{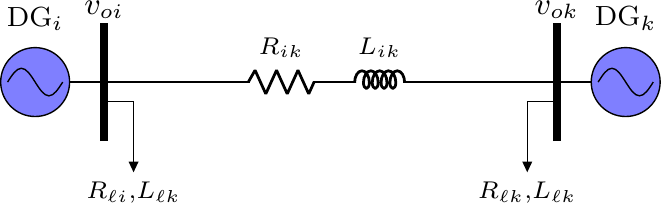}
	  		\vspace{-2mm}
	  		\caption{Two-inverter based autonomous microgrid}
	  		\rmspace
	  		\label{fig1}
	  	\end{center}
	  \end{figure}

\begin{figure}[h]
	  	\begin{center}
	  		\includegraphics[width=0.3\textwidth]{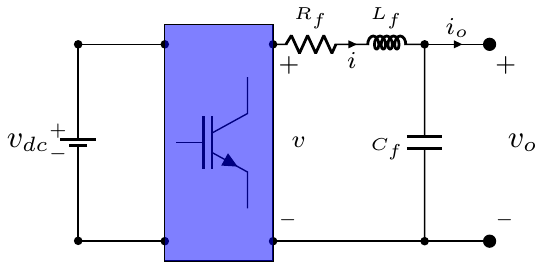}
	  		\vspace{-2mm}
	  		\caption{Schematic diagram of an AC inverter}
	  		\vspace{-3mm}
	  		\label{fig2}
	  	\end{center}
	  \end{figure}

For the \icl{comparison analysis} we consider the smallest unit and basic building block  of larger microgrids,  which is a typical two-inverter setup depicted by Fig. \ref{fig1}.
Our justification for using this setup is based on the fact that reduced-order models that do not correctly predict stability  for this basic configuration would also not be expected to perform  better in the  stability assessment of larger microgrids.
	 \icl{In this test} system (Fig. \ref{fig1}) each inverter is equipped with \icl{active} power based frequency droop \icl{control,} and reactive power based voltage \icl{droop control} \icl{(a widely used control policy as in e.g. \cite{pogaku2007, chandorkar1993, Guerrero2015, Guerrero2014,vorobev2017,guerrerodynamic2014, schiffer2014, SCHIFFER2016135, Bidram2013, Bidram2014, sao2005,ysunGuerrero2017,coelho2002, coelho2000})}.
	The inverters are connected to buses $i, k$ respectively via a balanced resistive-inductive line defined by the line parameters $R_{ik}, L_{ik} \in\mathbb{R}_{>0}$. Two local
	 resistive-inductive (RL) loads $R_{\ell i}, L_{\ell i} \in\mathbb{R}_{>0}$ and $R_{\ell k}, L_{\ell k} \in\mathbb{R}_{>0}$	  are connected to buses $i, k$ respectively. Fig. \ref{fig2} shows that the output of an inverter is fitted with an $LC$ filter with parameters $R_f,L_f,C_f\in\mathbb{R}_{>0}$. This attenuates the harmonics that are introduced into the output voltage by the pulse-width modulation scheme. The voltage across the capacitor $v_{oi}$ (respectively $v_{ok}$) determines the voltage of bus $i$ (respectively $k$).

\icl{The analysis in the rest of the paper is 
as follows. In section \ref{models} we provide a detailed description of the various models and derive those from first principles. 

In section \ref{accuracy} we evaluate the accuracy with which the reduced order models predict stability and instability. In particular, we find the range of values of the frequency and voltage droop gains for which each model is stable.  
These are then compared to those of the detailed average model, for different line $R/X$ ratios. Dynamic simulations of more advanced switching models are also presented for further validation.  Our findings demonstrate significant discrepancies  in the stability predictions of the reduced-order models and detailed average model, and further discussion on this is included in section~\ref{disc}.}

	\section{Autonomous Microgrid Dynamic Models}
\label{models}
In this section we \icl{describe in detail the models outlined in the previous section.}
 \icl{In particular, in order to provide intuition on their relative merits, these are rigourously derived from first principles stating throughout the approximations made for the simpler models to be deduced.}
	\subsection{Detailed Average Model}
\label{detailed}
The detailed average model we consider in our study consists of all the internal states of an inverter, and \icl{the} interconnection dynamics \cite{pogaku2007}. 
The inverter model uses the following standard assumptions: the DC-side (Fig. \ref{fig2}) is  equipped with sufficient energy reserves; and a sufficiently high switching frequency allows \icl{to neglect} 
the switching process.
These allow the inverter dynamics to be described by \icl{a continuous} average model which for convenience is formulated in the local $dq$ reference frame (Definition \ref{def3a}). In particular, the LC filter (Fig. \ref{fig2}) equations for inverter $i$ (respectively $k$) are obtained by applying $T(\theta_i)$ with $\dot\theta_i=\omega_i$ to the fundamental inductor and capacitor equations (e.g., \cite{pogaku2007, wang2016harmonic}) to give:
\begin{IEEEeqnarray}{rCl}
	\label{LC1}
	L_{fi}\dot i_{dqi}&=&\omega_iL_{fi}{J}i_{dqi}-R_{fi}i_{dqi}+v_{dqi}-v_{o{dqi}}\\
	\label{LC2}
	C_{fi} \dot v_{o{dqi}}&=&\omega_i	C_{fi}{J}v_{o{dqi}}+i_{dqi}-i_{o{dqi}}
\end{IEEEeqnarray}
where signals $i_{dqi}=[i_{di}~i_{qi}]^{\top}$,  $i_{odqi}=[i_{odi}~i_{oqi}]^{\top}$, $v_{dqi}=[v_{di}~v_{qi}]^{\top}$, $v_{odqi}=[v_{odi}~v_{oqi}]^{\top}$ take values in $\mathbb{R}^2$ and are the inverter $dq$ currents and voltages of the form described by \eqref{eqn5a}; $\omega_i$ is the inverter local frequency and takes value in $\mathbb{R}_{>0}$.
The inverter frequency and voltage set points are specified via the frequency and voltage droop controllers, which for inverter $i$ (respectively $k$) are given by \cite{pogaku2007, schiffer2014}:
\begin{subequations}\label{dt}
	\begin{align}\small
	\label{dt1}
	\tau\dot \omega_i = & -\omega_i + \omega_n - k_{pi}P_i  \\
	\label{dt2}
	\tau\dot V_i  = & -V_i + V_n - k_{qi}Q_i
	\end{align}
\end{subequations}
where $V_i$ is the inverter voltage amplitude and takes value in $\mathbb{R}_{>0}$;   $\tau, \omega_n, V_n, k_{pi}, k_{qi} \in \mathbb{R}_{>0}$ are respectively the filtering time constant, nominal frequency, nominal  voltage, frequency and voltage droop gains; $P_i$ and $Q_i$ are the respective measured active and reactive power.
The inverter output voltage $v_{odqi}$  is regulated to the voltage value set by \eqref{dt2} via the control input $v_{dqi}$. This is obtained from \icl{the} 
(outer) voltage and (inner) current controllers
which are respectively described for inverter $i$ (respectively $k$) as follows \cite{pogaku2007,Rocabert2012,teodo2004, wang2016harmonic}:
\begin{equation}
	\label{innerv}
	\begin{split}
	\dot \phi_{dqi}=&\mathbf{e}V_i-v_{odqi}\\
	\dot \gamma_{dqi}=&i_{dqi}^{ref}-i_{dqi}\\
	i_{dqi}^{ref}=&K_{PVi}(\mathbf{e}V_i-v_{odqi})+K_{IVi}\phi_{dqi}\\
	v_{dqi}=&K_{PCi}(i_{dqi}^{ref}-i_{dqi})+K_{ICi}\gamma_{dqi}+\omega_i L_{fi}Ji_{dqi}
	\end{split}
\end{equation}
where $\phi_{dqi}=[\phi_{di}~\phi_{qi}]^{\top}$, $\gamma_{dqi}=[\gamma_{dqi}~\gamma_{dqi}]^{\top}$ are the states of the voltage and current controllers respectively and take values in $\mathbb{R}^2$;  $K_{PVi}, K_{IVi}\in\mathbb{R}_{>0}$ are the proportional and integral gain of the voltage controller respectively; $K_{PCi}, K_{ICi}\in\mathbb{R}_{>0}$ are the proportional and integral gain of the current controller respectively; $i_{dqi}^{ref}=[i_{di}^{ref}~i_{qi}^{ref}]^{\top}$ is the reference current that the current controller tracks and \icl{takes} values in $\mathbb{R}^2$.

To facilitate the interconnection of the inverters, it is convenient to transform the port variables to the $DQ$ frame. The $DQ$ frame rotates with the synchronous frequency $\omega_{0}$ to which the inverter frequencies $\omega_i$ $(\omega_k)$ converge at steady state. \icl{Recalling}
Definition \ref{DQ0}, the RL line (Fig. \ref{fig1}) dynamics are derived by applying $T(\theta_0)$ to the fundamental inductor equation (e.g., \cite{pogaku2007}):
\begin{equation}
\begin{split}
	\label{dt3}
	L_{ik} \dot I_{DQ,ik} =& -R_{ik} I_{DQ,ik} + \omega_0 L_{ik}J I_{DQ,ik} \\& + \mathcal{T}(\delta_i)v_{odqi} - \mathcal{T}(\delta_k)v_{odqk}
	\end{split}\end{equation}
where $I_{DQ, ik}=[I_{D, ik}~I_{Q, ik}]^{\top}$ is the line current and takes values in $\mathbb{R}^2$, and
\begin{equation} \small
	\label{angle}
	\dot \delta_i  =  \omega_i - \omega_0.
\end{equation}	
where $\delta_i\in\mathbb{S}$ is the angle between the two reference frames.

The load connected to bus $i$ (respectively $k$) is considered as an RL load  and described on the synchronous reference frame \cite{pogaku2007}:
\begin{equation}\label{loaddq}
L_{\ell i} \dot I_{DQ,\ell i} = (-R_{\ell i}\mathbf{I}_2 + \omega_0 JL_{\ell i}) I_{DQ,\ell i} + \mathcal{T}(\delta_i)v_{odqi}
\end{equation}
where $I_{DQ,\ell i}=[I_{D,\ell i}~I_{Q,\ell i}]^{\top}$ is the load current and takes values in $\mathbb{R}^2$; $R_{\ell i}, L_{\ell i}\in\mathbb{R}_{>0}$ are the load resistance and inductance respectively.

In order to obtain $P_i$ and $Q_i$ in \eqref{dt}, the active and reactive power drawn by the load and that  transferred over the line  are required.
The  power drawn by the load at bus $i$ (respectively $k$) is given by
\begin{equation}\label{loadpq}
\begin{split}
P_{\ell i}=&I_{DQ,\ell i}^{\top}\mathcal{T}(\delta_i)v_{odqi}\\
Q_{\ell i} =& I_{DQ,\ell i}^{\top}J\mathcal{T}(\delta_i)v_{odqi}
\end{split}
\end{equation}
where $P_{\ell i}$, $Q_{\ell i}$ are the active and reactive power respectively.
The active and reactive power flows over the line from bus $i$ to bus $k$ are given as
\begin{equation}
	\begin{split}
	P_{ik}=&I_{DQ,ik}^{\top}\mathcal{T}(\delta_i)v_{odqi}\\
	Q_{ik} =& I_{DQ,ik}^{\top}J\mathcal{T}(\delta_i)v_{odqi}.
	\end{split}
\end{equation}
The power injected by inverter $i$ is then computed as:
\begin{equation*}
	\begin{split}
	P_i=&P_{\ell i}+P_{ik}=(I_{DQ,ik}+I_{DQ,\ell i})^{\top}\mathcal{T}(\delta_i)v_{odqi} \\
	Q_i=&Q_{\ell i}+Q_{ik}=(I_{DQ,ik}+I_{DQ,\ell i})^{\top}J\mathcal{T}(\delta_i)v_{odqi}
	\end{split}
\end{equation*}
and can also be written as
\begin{equation}	\label{pq}
	\begin{split}
	P_i =& i_{odqi}^{\top}v_{odqi}=(I_{DQ,ik}+I_{DQ,\ell i})^{\top}\mathcal{T}(\delta_i)v_{odqi} \\
	Q_i =& i_{odqi}^{\top}Jv_{odqi}=(I_{DQ,ik}+I_{DQ,\ell i})^{\top}J\mathcal{T}(\delta_i)v_{odqi}
	\end{split}
\end{equation}
Recall that ${J}\mathcal{T}(\delta_i)=\mathcal{T}(\delta_i){J}$ (Remark \ref{def8}), thus it clearly follows from \eqref{pq} that
\begin{equation}
\label{iodq}
i_{odqi}=\mathcal{T}^{-1}(\delta_i)(I_{DQ,ik}+I_{DQ,\ell i}).
\end{equation}
Hence the detailed average model is described by \eqref{LC1}--\eqref{angle}, \eqref{iodq}, which is used together with the load equation \eqref{loaddq}, and injected power \eqref{pq}.

	\subsection{Electromagnetic 5th-Order Model}
	\label{5th}
	\icl{In the} detailed average model presented in \ref{detailed} \icl{the} voltage across the capacitor (Fig. \ref{fig2}) tracks a reference via the control action of the outer and inner controllers.
\icl{The electromagnetic (EM) 5th-order model is formulated by making the assumption} that the outer and inner controllers are tuned to be much faster than the droop \icl{controllers 
such that} a fast tracking response is obtained.  
By this reasonable assumption (similarly applied in \cite{Guerrero2015, Guerrero2014, guerrerodynamic2014}), dynamics \eqref{LC1}-\eqref{LC2}, \eqref{innerv}
are neglected, and the inverter is considered as a voltage source which produces at its output (Fig. \ref{fig1}) a symmetric three-phase AC voltage
	$v_{oi}:\mathbb{R}_{>0} \to \mathbb{R}^3$ of the form
	\begin{equation}\small
	\label{vabc}
	v_{{oi}}={V_i\begin{bmatrix} \sin(\theta_i) \\ \sin(\theta_i-\frac{2\pi}{3}) \\\sin(\theta_i+\frac{2\pi}{3})) \end{bmatrix}}
	\end{equation}
	where $\theta_i$ is the angle defined as $\dot\theta_i=\omega_i$, and $V_i$, $\omega$ are respectively the voltage amplitude and frequency and take values in $ \mathbb{R}_{>0}$. Due to this symmetry, $\mathds{1}^{\top}_3v_{oi}=0$.
	Thus, according to \eqref{eqn5a}, \eqref{vabc} becomes
	\begin{equation}
	\label{vdq}
		v_{o{dq_i}}=T(\theta_i)v_{oi}=\mathbf{e}V_i.
	\end{equation}
	For interconnection with the network, it is convenient to map the inverter output voltage to the synchronous rotating reference frame (Remark \ref{def8}) as follows:
	\begin{equation}\label{vodq}
	\mathcal{T}(\delta_i)v_{o{dq_i}}=\mathcal{T}(\delta_i)\mathbf{e}V_i= \begin{bmatrix}\cos(\delta_i) \\\sin(\delta_i)	\end{bmatrix}V_i
	\end{equation}
\icl{Applying the model simplifications described above}, 
the electromagnetic (EM) 5th-order model with three states relating to the $i^{th}$ (respectively $k^{th}$) inverter (i.e. $\delta_i, \omega_i, V_i$) and two states for the lines ($DQ$ components of current $I_{DQ,ik}$) \cite{vorobev2018, Guerrero2015} is:
\begin{subequations} \small
\label{5thh}
	\begin{align}
	\label{5th1}
	\dot \delta_i  = & \omega_i - \omega_0  \\
	\label{5th2}
	\tau\dot \omega_i = & -\omega_i + \omega_n - k_{pi}P_i  \\
	\label{5th3}
	\tau\dot V_i  = & -V_i + V_n - k_{qi}Q_i  \\
	\label{5th4}
	L_{ik} \dot I_{D,ik} = & -R_{ik} I_{D,ik} + \omega_0 L_{ik} I_{Q,ik} + V_i\cos\delta_i - V_k\cos\delta_k \\
	\label{5th5}
	L_{ik} \dot I_{Q,ik} = & -R_{ik} I_{Q,ik} - \omega_0 L_{ik} I_{D,ik} + V_i\sin\delta_i - V_k\sin\delta_k
	\end{align}
\end{subequations}
Considering \eqref{vodq} the load model \eqref{loaddq} is rewritten as
\begin{equation}\label{loaddq2}
L_{\ell i} \dot I_{DQ,\ell i} = (-R_{\ell i}\mathbf{I}_2 + \omega_0 JL_{\ell i}) I_{DQ,\ell i} + \mathcal{T}(\delta_i)\mathbf{e}V_i.
\end{equation}
The injected power $P_i$, $Q_i$ in \eqref{5th2}--\eqref{5th3} are computed as \mbox{$P_i=P_{\ell i}+P_{ik}$} and \mbox{$Q_i=Q_{\ell i}+Q_{ik}$}.
Considering \eqref{vodq} \icl{we} rewrite \eqref{loadpq} to obtain the active and reactive power drawn by the load at bus $i$ (respectively $k$) as follows:
\begin{equation}\small
\label{pqload}
\begin{split}
P_{\ell i} &= I^{\top}_{DQ,\ell i}\mathcal{T}(\delta_i)\mathbf{e}V_i \\
Q_{\ell i} &= I^{\top}_{DQ,\ell i}J\mathcal{T}(\delta_i)\mathbf{e}V_i.
\end{split}
\end{equation}	
	Note that \icl{various} 
reduced-order models \icl{can differ based on} 
the approximations asserted on the line dynamics which also affects how the power exchange $P_{ik}, Q_{ik}$ are computed. 	For the EM 5th-order model which considers the lines dynamics \eqref{5th4}--\eqref{5th5}, 
and the corresponding  power flows over the line from bus $i$ to bus $k$ is given by:
	\begin{equation} \small
		\label{pik}
		\begin{split}
		P_{ik} &= I^{\top}_{DQ,ik}\mathcal{T}(\delta_i)\mathbf{e}V_i\\
		Q_{ik} &= I^{\top}_{DQ,ik}J\mathcal{T}(\delta_i)\mathbf{e}V_i
		\end{split}
	\end{equation}
	The injected power  $P_i$, $Q_i$ are
	\begin{equation}\label{pq5}
	\begin{split}
	P_i=&P_{\ell i}+P_{ik}=(I_{DQ,ik}+I_{DQ,\ell i})^{\top}\mathcal{T}(\delta_i)\mathbf{e}V_i\\
	Q_i=&Q_{\ell i}+Q_{ik}=(I_{DQ,ik}+I_{DQ,\ell i})^{\top}J\mathcal{T}(\delta_i)\mathbf{e}V_i.
	\end{split}
	\end{equation}
 \icl{The 5th-order model of an inverter is therefore given by \eqref{5thh}, which} is used together with \eqref{loaddq2}, \eqref{pq5}.

	\subsection{Conventional 3rd-Order Model}
	\label{3th}

\icl{The conventional 3rd-order model uses the same assumptions on the inverter dynamics as the 5th order model (section \ref{5th}), and hence the inverter dynamics are described by \eqref{5th1}--\eqref{5th3}.}
The 
\icl{main difference} 
between the 3rd-order and \icl{the} 5th-order model comes from the approximations asserted on the line dynamics,
	which also affects how the power exchange $P_{ik}, Q_{ik}$ are calculated.
The conventional 3rd-order model uses a traditional power system assumption of a distinct timescale separation which allows the line dynamics to be neglected and the lines modelled statically. In particular, the conventional 3rd-order model uses the traditional quasi-stationary approximation	
	(also referred to as zero-order approximation model) \cite{kundur94, machowski08, padiyar1996, anderson02, glover11, schiffer2014}.
	This involves evaluating  \eqref{5th4}--\eqref{5th5} at equilibrium and corresponds to setting the derivative terms $  \dot I_{D,ik}, \dot I_{Q,ik}$ to zero as follows: 	
	\begin{subequations} \label{45a}
		\begin{align}
		\small
		\label{5th4a}
	0 = & -R_{ik} I_{D,ik} + \omega_0 L_{ik} I_{Q,ik} + V_i\cos\delta_i - V_k\cos\delta_k \\
		\label{5th5a}
	0 = & -R_{ik} I_{Q,ik} - \omega_0 L_{ik} I_{D,ik} + V_i\sin\delta_i - V_k\sin\delta_k.
		\end{align}
	\end{subequations}
	For simplicity \icl{in the} analysis  \eqref{45a} is expressed as a phasor. To do this  rewrite \eqref{45a} by multiplying \eqref{5th5a} with the \yo{complex number} $j$ as follows
	\begin{subequations} \label{45b}
		\begin{align}
		\small
		\label{5th4b}
		0 = & -R_{ik} I_{D,ik} + \omega_0 L_{ik} I_{Q,ik} + V_i\cos\delta_i - V_k\cos\delta_k \\
		\label{5th5b}
		0 = & -jR_{ik} I_{Q,ik} - j\omega_0 L_{ik} I_{D,ik} + jV_i\sin\delta_i - jV_k\sin\delta_k.
		\end{align}
	\end{subequations}
	The summation of \eqref{5th4b} and \eqref{5th5b} gives:
	\begin{equation}\small
	\label{io}
		I^0_{ik}=(R_{ik}+j\omega_0L_{ik})^{-1}(V_i e^{j\delta_i}-V_ke^{j\delta_k})
	\end{equation}
	where $I^0_{ik}=I_{D,ik}+jI_{Q,ik}$ is the current phasor and the \icl{superscript $0$ denotes} that it is calculated at zero-order approximation.
	The corresponding zero-order approximation of the active and reactive power exchange is obtained via the relationship
	\begin{equation*}\small
		P^0_{ik}+jQ^0_{ik}=(I^0_{ik})^*V_i e^{j\delta_i}
	\end{equation*}
and this simplifies to:
	\begin{equation}\small
		\label{pq0}
		\begin{split}
		P^0_{ik} &=  G_{ik}V^2_i - G_{ik}V_iV_k\cos(\delta_i-\delta_k)   +B_{ik}V_iV_k\sin(\delta_i-\delta_k)  \\
		Q^0_{ik} &=  B_{ik}V^2_i - B_{ik}V_iV_k\cos(\delta_i-\delta_k)
		- G_{ik}V_iV_k\sin(\delta_i-\delta_k)
		\end{split}
	\end{equation}
	where
	\begin{equation*}\small
		G_{ik}=\frac{R_{ik}}{R^2_{ik}+X^2_{ik}}, ~~ B_{ik}=\frac{X_{ik}}{R^2_{ik}+X^2_{ik}}, ~X_{ik}=\omega_0L_{ik},
	\end{equation*}
	$G_{ik}, B_{ik}, X_{ik} \in \mathbb{R}_{>0}$ are the conductance, susceptance and reactance respectively.

	The injected power $P_i$, $Q_i$  are computed as
	\begin{equation}\label{pq11}
	\begin{split}
	P_i=&P_{\ell i}+P_{ik}^0\\
	Q_i=&Q_{\ell i}+Q_{ik}^0
	\end{split}
	\end{equation}
	where $P_{\ell i}$, $Q_{\ell i}$ are given by \eqref{pqload}, and $P_{ik}^0$, $Q_{ik}^0$ obtained from \eqref{pq0}.\newline
	Hence  \icl{the conventional 3rd-order model is described by 
\eqref{5th1}--\eqref{5th3}, \eqref{pq0}},
which is used together with \eqref{loaddq2}--\eqref{pqload}, \eqref{pq11}.
	
	\subsection{High-Fidelity 3rd-Order Model}
	\label{hf}
	The high-fidelity 3rd-order model was proposed in \cite{vorobev2018} to \icl{improve the conventional} 
3rd-order model. The authors in \cite{vorobev2018} raised the concern that the inherent low inertia of inverters may not allow for straight-forward assumption of neglecting the line dynamics.
	The fast dynamics of lines can influence the slow ones despite their short timescale.  This is in contrast to
\icl{conventional power systems} where a distinct timescale separation exists.
	
	\icl{Below we describe the high-fidelity 3rd-order model following its derivation in \cite{vorobev2018}; 
in the description we relax the assumption that}
an inverter 
\icl{is} connected to a fixed voltage bus.
\icl{The high-fidelity 3rd-order model uses the same assumptions on the inverter dynamics \icl{as the} 5th order model  (section \ref{5th}), and hence inverter dynamics are described by \eqref{5th1}--\eqref{5th3}.}
	The distinction between the  high-fidelity 3rd-order model and the other two reduced-order models comes from the approximations asserted on the line dynamics, which also affects how the power exchange $P_{ik}, Q_{ik}$ are computed.
	Instead of setting the derivative terms  $ \dot I_{D,ik}, \dot I_{Q,ik}$  in \eqref{5th4} and \eqref{5th5} to zero,
	\yo{we take  \icl{Laplace transforms} which leads to the following\footnote{We slightly abuse notation by using the same \icl{symbol for a time domain variable and its Laplace transform.}}:}
	\begin{subequations} \label{45c}
		\begin{align}
		\small
		\label{5th4c}
		sL_{ik} I_{D,ik} = & -R_{ik} I_{D,ik} + \omega_0 L_{ik} I_{Q,ik} + V_i\cos\delta_i \notag\\&- V_k\cos\delta_k \\
		\label{5th5c}
	sL_{ik} I_{Q,ik} = & -R_{ik} I_{Q,ik} - \omega_0 L_{ik} I_{D,ik} + V_i\sin\delta_i \notag \\&- V_k\sin\delta_k
		\end{align}
	\end{subequations}
	For convenience \icl{in the} analysis \icl{we} rewrite \eqref{45c} by multiplying \eqref{5th5c} with the \yo{complex number}  $j$ to obtain
	\begin{subequations} \label{45d}
		\begin{align}
		\small
		\label{5th4d}
		sL_{ik} I_{D,ik} = & -R_{ik} I_{D,ik} + \omega_0 L_{ik} I_{Q,ik} + V_i\cos\delta_i \notag\\&- V_k\cos\delta_k \\
		\label{5th5d}
		jsL_{ik} I_{Q,ik} = & -jR_{ik} I_{Q,ik} - j\omega_0 L_{ik} I_{D,ik} +j V_i\sin\delta_i \notag\\&- jV_k\sin\delta_k
		\end{align}
	\end{subequations}
	 The summation of \eqref{5th4d} and \eqref{5th5d} gives
	\begin{equation*}\small
	\begin{split}
	I^1_{ik}=&		\frac{V_i e^{j\delta_i}-V_ke^{j\delta_k}}{R_{ik}+j\omega_0L_{ik}+sL_{ik}}
	\\=& \frac{(R_{ik}+j\omega_0L_{ik})^{-1}(V_i e^{j\delta_i}-V_ke^{j\delta_k})}{1+\tfrac{sL_{ik}}{R_{ik}+j\omega_0L_{ik}}}
	\end{split}
	\end{equation*}
	and this can be simplified to
	\begin{equation}
		\label{fI1}
		I^1_{ik}=		\frac{I^0_{ik}}{1+\tfrac{sL_{ik}}{R_{ik}+j\omega_0L_{ik}}}
	\end{equation}
	where   $I^1_{ik}=I_{D,ik}+jI_{Q,ik}$ is the current phasor.
	For the derivation of the equivalent reduced-order model that captures the dynamics of the slow modes, a reasonable assumption is that \mbox{$|\tfrac{L_{ik}}{R_{ik}+j\omega_0L_{ik}}|$} is sufficiently small compared to the
	 electromagnetic time $\frac{L_{ik}}{R_{ik}}$ such that
\eqref{fI1} is represented by the first-order approximation of the Taylor expansion as follows:
	\begin{equation}\small
	\label{fI2}
		I^1_{ik}\approx I^0_{ik}- \frac{L_{ik}}{R_{ik}+j\omega_0 L_{ik}} sI^0_{ik}.
	\end{equation}
	Returning to the time domain \eqref{fI2} can be rewritten as
	\begin{equation}\small
	\label{fI3}
		I^1_{ik}\approx I^0_{ik}- \frac{L_{ik}}{R_{ik}+j\omega_0 L_{ik}} \dot I^0_{ik}
	\end{equation}
	where $I^0_{ik}$ is as in \eqref{io}, and the term $\dot I^0_{ik}$ reads
	\begin{equation*}\small
		\dot I^0_{ik}=\dot V_ie^{j\delta_i}+jV_i\dot \delta_i e^{j\delta_i}-\dot V_k e^{j\delta_k}-jV_k\dot \delta_k e^{j\delta_k}.
	\end{equation*}
	The corresponding first-order approximation of the active and reactive power is obtained from
	\begin{equation}\small
	P^1_{ik}+jQ^1_{ik}\approx(I^1_{ik})^*V_i e^{j\delta_i}
	\end{equation}
	which simplifies to:
	\begin{equation}\small
	\label{pq1a}
	\begin{split}
		P^1_{ik} \approx ~&  P^0_{ik} - G^\prime_{ik}V_i\dot V_i + G^\prime_{ik}V_i\dot V_k \cos(\delta_i-\delta_k) \\&  + G^\prime_{ik}V_i V_k \dot \delta_k \sin(\delta_i-\delta_k) - B^\prime V^2_i \dot \delta_i \\& -B^\prime_{ik}V_i\dot V_k \sin(\delta_i-\delta_k) + B^\prime_{ik}V_i V_k \dot \delta_k \cos(\delta_i-\delta_k)     	\end{split}
\end{equation}
	\begin{equation}\small
	\label{pq1b}
		\begin{split}
		Q^1_{ik} \approx~ &   Q^0_{ik} + G^\prime_{ik}V^2_i\dot \delta_i + G^\prime_{ik}V_i\dot V_k \sin(\delta_i-\delta_k)  \\&  - G^\prime_{ik}V_i V_k \dot \delta_k \cos(\delta_i-\delta_k) - B^\prime V_i \dot V_i \\& + B^\prime_{ik}V_i\dot V_k \cos(\delta_i-\delta_k) + B^\prime_{ik}V_i V_k \dot \delta_k \sin(\delta_i-\delta_k)     	\end{split}
\end{equation}
	where
	\begin{equation*}
		G^\prime_{ik}=\frac{(R^2_{ik}-X^2_{ik})L_{ik}}{(R^2_{ik}+X^2_{ik})^2},
		~~ B^\prime_{ik}=\frac{2R_{ik}X_{ik}L_{ik}}{(R^2_{ik}+X^2_{ik})^2},
	\end{equation*}
	$G^\prime_{ik}\in\mathbb{R}$, $B^\prime_{ik}\in\mathbb{R}_{>0}$ can be referred to as the subsynchronous conductance and susceptance \cite{vorobev2018}. Typical values of $G^\prime_{ik}, B^\prime_{ik}$ are small, but could play an important role during the transient \icl{response} of the voltages and angles. Note that $P^0_{ik}$, $Q^0_{ik}$ given by \eqref{pq0} in the conventional 3rd-order model can be recovered if $G^\prime_{ik}, B^\prime_{ik}$ are set to zero.
	
	The injected power $P_i$, $Q_i$  are computed as
	\begin{equation}\label{pq1}
	\begin{split}
	P_i= & P_{\ell i}+P_{ik}^1\\
	Q_i= & Q_{\ell i}+Q_{ik}^1
	\end{split}
	\end{equation}
	where $P_{\ell i}$, $Q_{\ell i}$ are given by \eqref{pqload}, and $P_{ik}^1$, $Q_{ik}^1$ obtained from \eqref{pq1a}, \eqref{pq1b} respectively.\newline
	Hence  dynamics \eqref{5th1}--\eqref{5th3}, \eqref{pq1a}--\eqref{pq1b} describe the  high-fidelity 3rd-order model and \icl{are} used together with \eqref{loaddq2}--\eqref{pqload}, \eqref{pq1}.

	\begin{table}[ht]
		\centering
		\vspace{-0mm}
		\caption{Microgrid parameters}\label{tb1}
		\scriptsize
		\begin{tabular}{c|c}
			Description & Value \\
			\hline
			Inverter $i$, $k$ & $R_{f}$ = 0.1 $\Omega$, $L_{f}$ = 5 $\text{mH}$, $C_{f}$= 50 $\mu$F,
			\\ &  $k_{p}$ = 6 $\times$ $10^{-5}$, $k_{q}$ = 1.5 $\times$ $10^{-4}$, $K_{PV}$ = 5, \\&$K_{IV}$ = 10, $K_{PC}$ = 5, $K_{IC}$ = 25, \\&$\omega_n$ = 2$\pi$(50) rad/s, $\tau$ = 31.8 ms, $V_n$ = 311 V.\\
			\hline
			R/X ratios & R/X $\gg$ 1 : $R_{ik}$ = 641 m$\Omega$, $L_{ik}$ = 0.26 mH
			 \\& R/X $\approx$ 1 : $R_{ik}$ = 195 m$\Omega$, $L_{ik}$ = 0.61 mH
			 \\& R/X $\ll$ 1 : $R_{ik}$ = 0.4 $\Omega$, $L_{ik}$ = 7 mH\\
			\hline
		Loads &  $R_{\ell i}$ = 20 $\Omega$, $R_{\ell k}$ = 40 $\Omega$, $L_{\ell i}$ = 15 mH, $L_{\ell k}$ = 40 mH
		\end{tabular}
		\label{tab:params}
	\end{table}

	\section{Small-Signal Models}
	\label{smallsignal}
	The \icl{stability analysis} 
in our comparative study in section \ref{accuracy} \icl{is based on a small signal analysis} 
of the detailed average model and reduced-order models.
 In order to obtain the small-signal models,  we  linearized  the nonlinear dynamics  (section \ref{models})  around an equilibrium point $x^s \in \{\delta^s_i, \omega_0, V^s_i, \phi^s_{dqi}, \gamma^s_{dqi}, i^s_{dqi},v^s_{odqi}, I^s_{DQ,ik},I^s_{DQ,\ell i}\}$.
 The equilibrium \icl{points are} obtained from \icl{simulations} of the nonlinear models (section \ref{models}) in three different regimes of line $R/X$ ratios for the setup in Fig. \ref{fig1}. The system parameters are given in Table \ref{tab:params}. The small-signal linearization of the detailed average model and reduced-order models can be found  \icl{ in the Appendix.}

	\section{Model Accuracy Assessment}
	\label{accuracy}	
	In this section, we assess the quality of  the  EM 5th-order, conventional 3rd-order, and high-fidelity 3rd-order  models via the comparison of their  stability \icl{properties to those} of the detailed average model \icl{(section \ref{detailed})}.

	\subsection{Small-Signal Stability Prediction}
	\label{ss}

\icl{For the comparison in this section we deduce the range of values of the frequency and voltage droop gains ($k_p$ , $k_q$ respectively) for which each of the models is stable. This is carried out by increasing the droop gains until instability occurs.}
\icl{We will refer to the range of values of the droop gains where stability is maintained as the {\em stability region} of the model.}
\icl{These regions are deduced for each model for three different values of the line $R/X$ ratios, and are compared to those of the detailed model, as a means of evaluating the accuracy of the reduced order models.}

\icl{Stability is determined via a small signal analysis and the corresponding eigenloci are also presented.}
	The $R/X$ ratio \icl{regimes} considered are $R/X\gg1$, $R/X\approx1$, and $R/X\ll1$ respectively, which \icl{are regimes} 
encountered in \icl{microgrids.}
	
	 It should be noted that the voltage and  current \icl{control} gains of the inverters in the detailed average model are 
tuned to obtain \icl{the best possible} performance. This is to allow a fair comparison of the detailed average model with the reduced-order \icl{models.}
Also note \icl{that the tests described in this section were repeated for various inverter/load/line parameters and similar conclusions were obtained. For convenience in the presentation the results are presented for typical equilibrium points with the parameters as described within the paper.}
	
\begin{figure}[ht!]
	\begin{center}
		\includegraphics[width=0.4\textwidth]{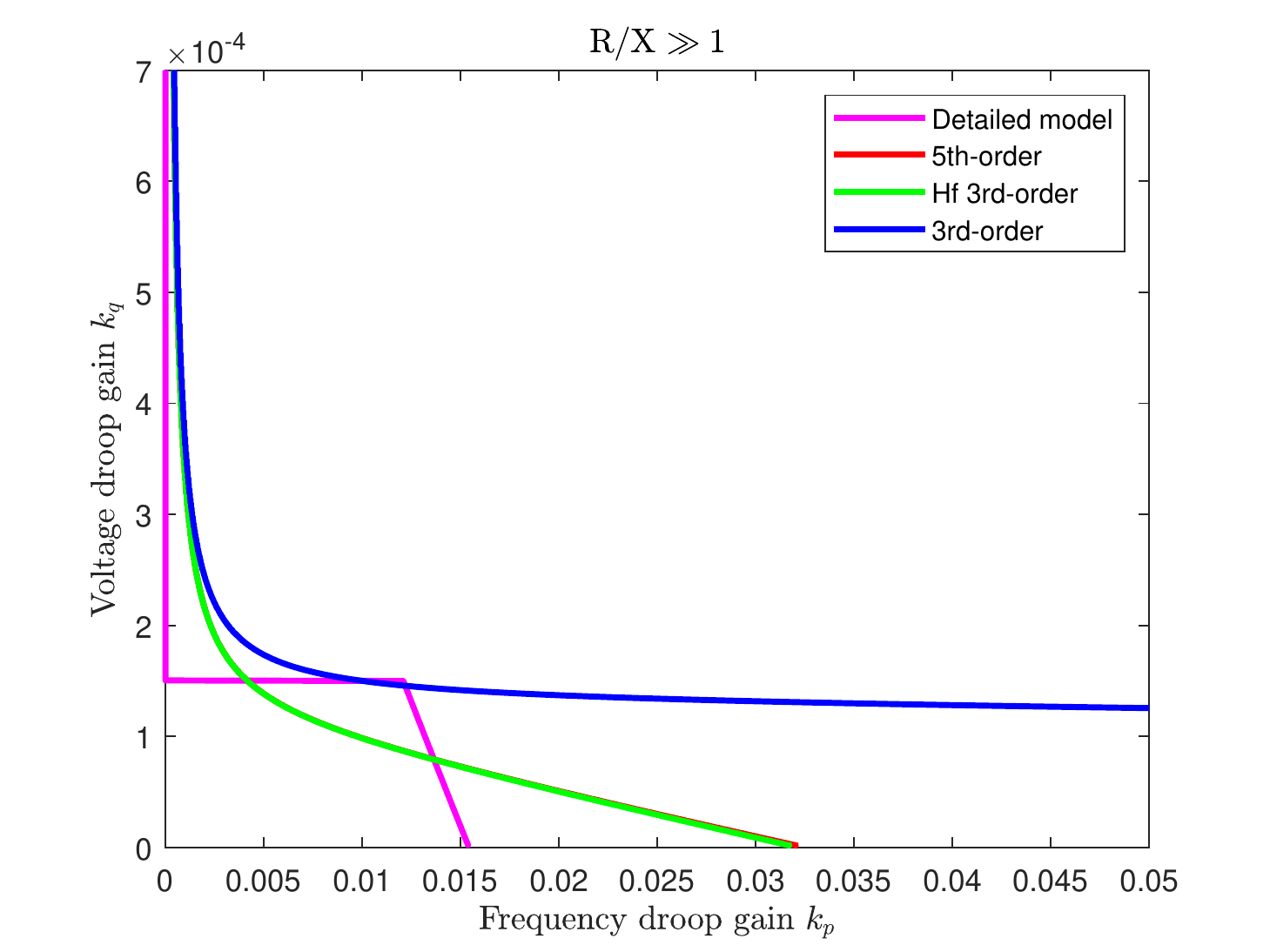}
		\vspace{-3mm}
		\caption{Comparison of stability regions given by the four small-signal  models for the case $R/X\gg1$. The region to the left of the corresponding curve are the values of the droop gains  for which the model is stable.}.
		\vspace{-2mm}
		\label{regioLV}
	\end{center}
\end{figure}

\begin{figure}[ht!]
	\begin{center}
		\includegraphics[width=0.4\textwidth]{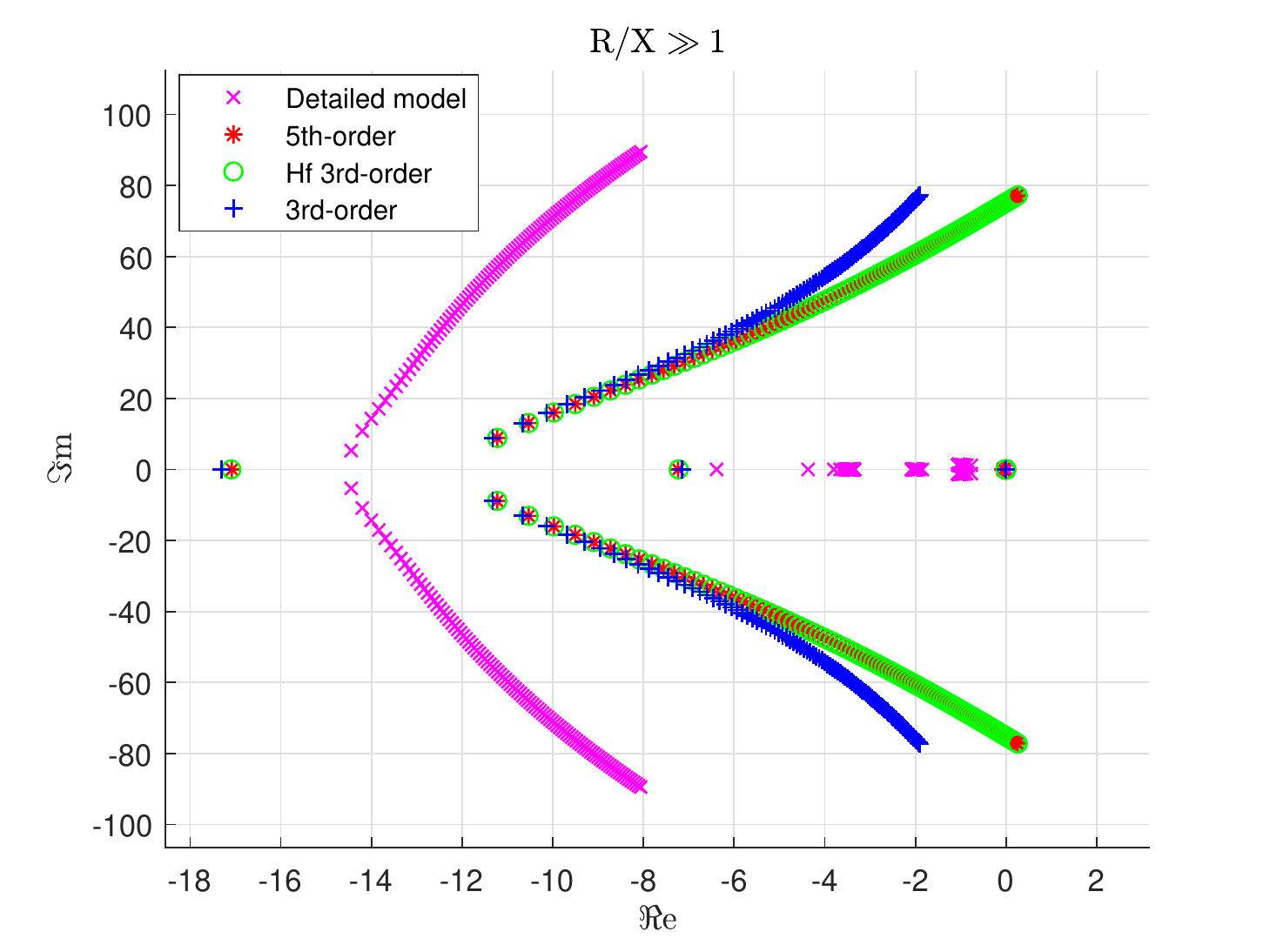}
		\vspace{-3mm}
		\caption{Eigenloci of the four  small-signal models for the case $R/X\gg1$ when \icl{$k_{pi}=k_{pk}$ take values in the range $[6\times10^{-5}, 4.4\times10^{-3}]$, and} \icl{$k_{qi}=k_{qk}=1.5\times10^{-4}$}.}
		\vspace{-2mm}
		\label{LV}
	\end{center}
\end{figure}

\begin{figure}[ht!]
	\begin{center}
		\includegraphics[width=0.4\textwidth]{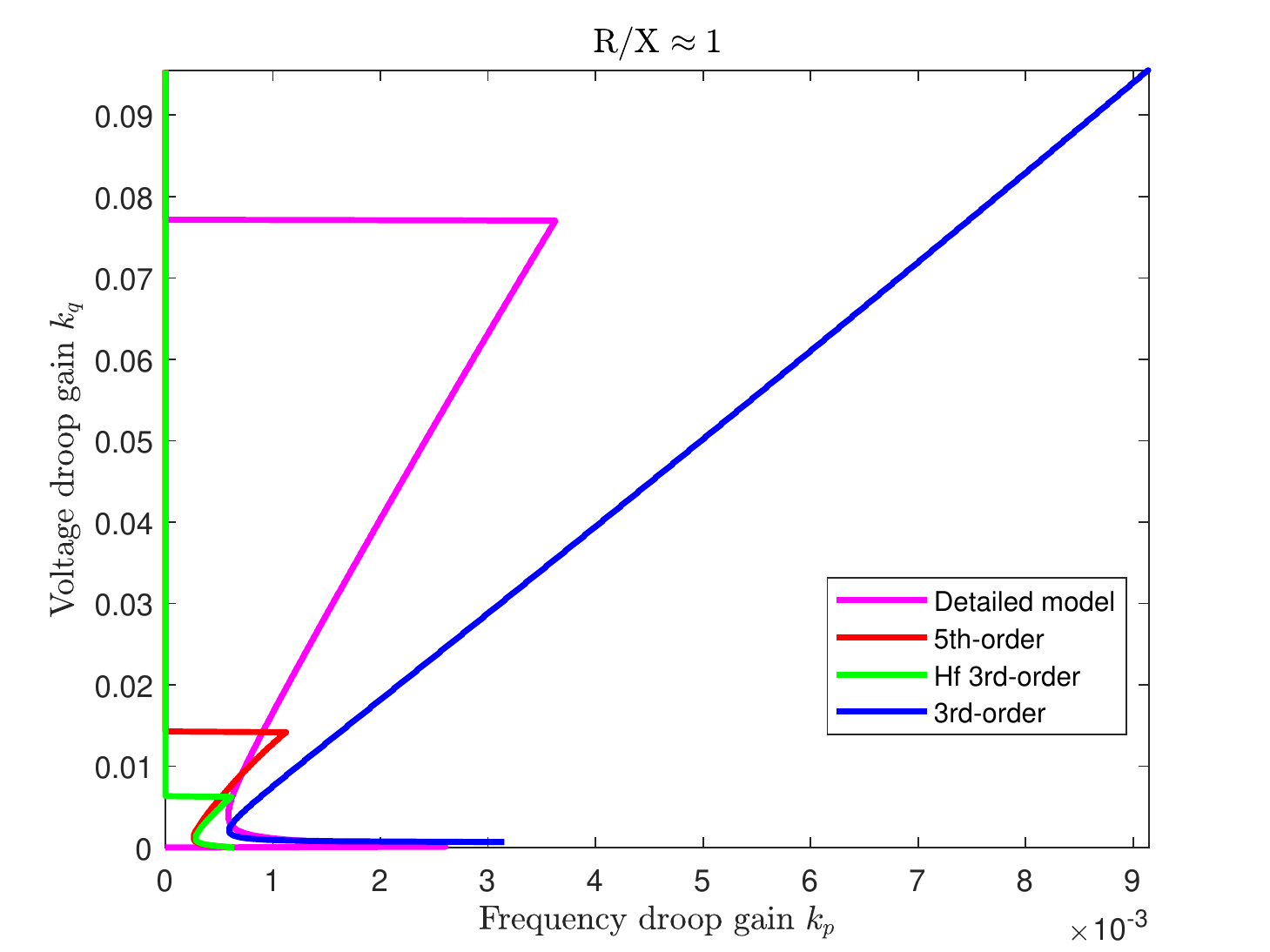}
		\vspace{-3mm}
		\caption{Comparison of stability regions given by the four small-signal models for the case $R/X\approx1$. The region to the left of the corresponding curve are the values of the droop gains  for which the model is stable.}
		\vspace{-3mm}
		\label{regionMV}
	\end{center}
\end{figure}

\begin{figure}[ht!]
	\begin{center}
		\includegraphics[width=0.4\textwidth]{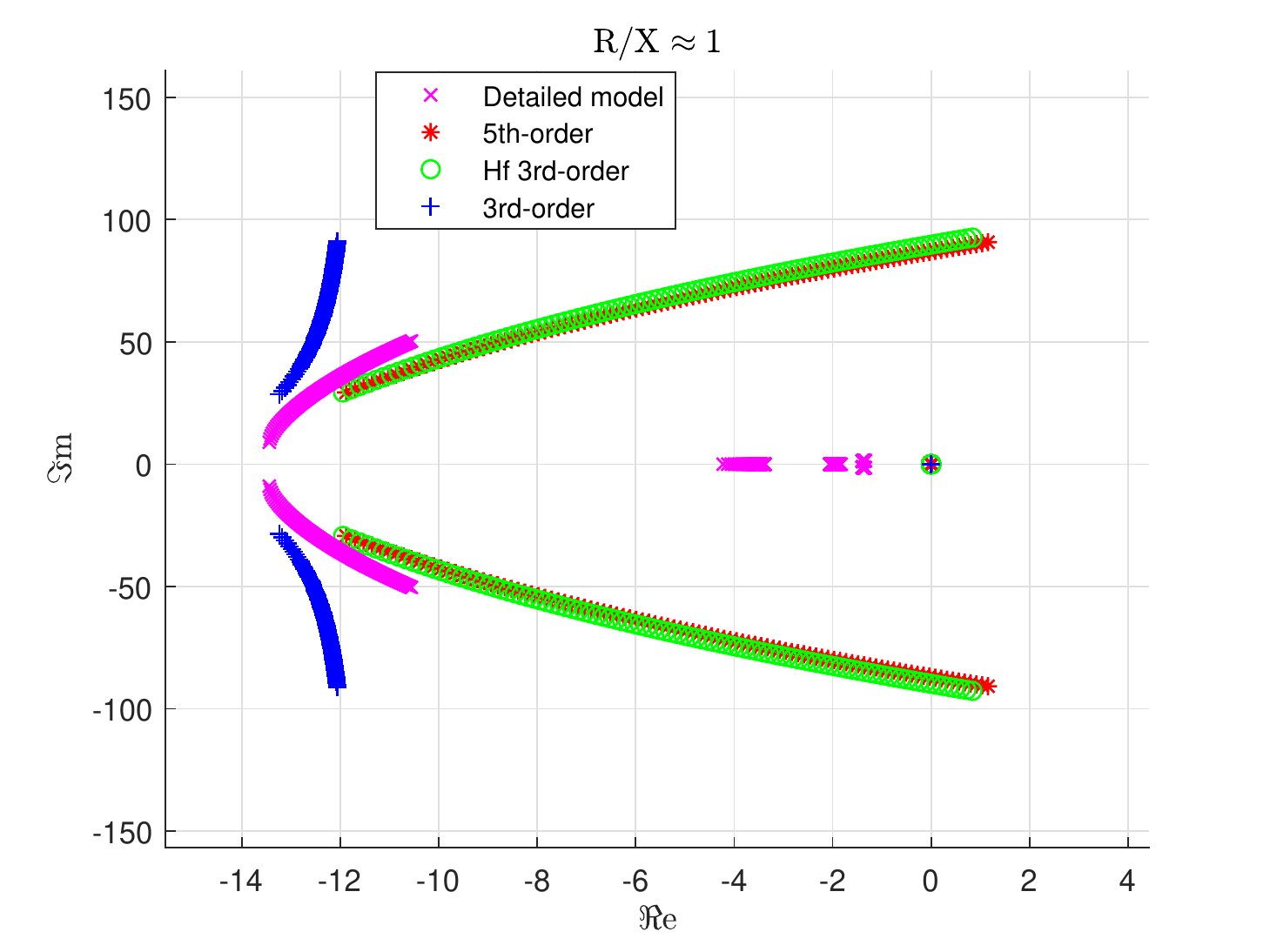}
		\vspace{-3mm}
		\caption{Eigenloci of the four small-signal models for the case $R/X\approx1$ when
			\icl{$k_{pi}=k_{pk}$ take values in the range $[6\times10^{-5}, 5.3\times10^{-4}]$ and
				$k_{qi}=k_{qk}=1.5\times10^{-4}$}.}
		\vspace{-3mm}
		\label{MV}
	\end{center}
\end{figure}

	We now present the \icl{comparison results}. Figs. \ref{regioLV}, \ref{regionMV}, \ref{regionHV} show the stability region \icl{of each model in the three different $R/X$ regimes}. 
\icl{In particular, the} region to the left of the \icl{each of the curves} are the values of the droop gains  for which the \icl{corresponding} model is stable.  Note that in the legend of Figs. \ref{regioLV}--\ref{HV} \icl{detailed} model refers to the detailed average \icl{model;} 5th-order refers to the EM 5th-order model; 3rd-order refers to the conventional 3rd-order model; and Hf 3rd-order refers to the high-fidelity 3rd-order model.

\vspace{-0mm}

	\subsubsection{Line ratio $R/X\gg1$} \icl{To investigate this regime we considered}
line parameters \mbox{$R_{ik}=641~\text{m}\Omega$,} \mbox{$L_{ik}=0.26~\text{mH}$},  $R/X\approx7.8$ \cite{engler2005}.
 	\paragraph{Stability Region}
\icl{The stability region of the four models is presented in Fig. \ref{regioLV}.}
	It is evident that the 5th-order and high-fidelity 3rd-order \icl{models have the same stability region.} 
 	Their \icl{stability region} is greater than that of the detailed average model, but much \icl{smaller} than that of the conventional 3rd-order model.
 	Compared to the detailed average model, all the reduced-order models give erroneous stability results \icl{for large} 
 frequency and voltage droop gains. The conventional 3rd-order model, \icl{in particular}, has \icl{the biggest discrepancy thus demonstrating that 
 it is 
 inappropriate} for stability assessment when $R/X\gg1$.

 	\paragraph{Eigenloci}
 	Fig. \ref{LV} shows the eigenloci \icl{of each of the four models \icl{(i.e. the eigenvalues of state matrix $A$ in a state space representation)}. In particular, the eigenloci are plotted as the frequency droop coefficients $k_{pi}=k_{pk}$ vary in the interval  $[6\times10^{-5}, 4.4\times10^{-3}]$. The voltage droop coefficients are set to $k_{qi}= k_{qk}=1.5\times10^{-4}$} .
 	It is clear that the eigenloci of the 5th-order and high-fidelity 3rd-order models are \icl{are very close}, which confirms their similar stability regions (Fig. \ref{regioLV}).
 \icl{It should also be noted that the location of the poles} of the reduced-order models \icl{does} not change much for \icl{large droop gains} compared to that of  the detailed average model. This causes the reduced models to give erroneous stability regions \icl{as was demonstrated in  Fig.~\ref{regioLV}.}

\subsubsection{Line ratio $R/X\approx 1$}
\icl{To investigate this regime we considered}
line parameters \mbox{$R_{ik}=195~\text{m}\Omega$,} \mbox{$L_{ik}=0.61~\text{mH}$}, \mbox{$R/X\approx1.0$} \cite{vorobev2018}.
	
\paragraph{Stability Region}
\icl{The stability region of each of the four models}
for the case $R/X\approx1$ is presented in Fig. \ref{regionMV}.
The \icl{stability region of the} high-fidelity 3rd-order  is within  that of the detailed average model. \icl{A similar behaviour is} observed with the 5th-order model, \icl{but this} admits \icl{a small} error for some values of the frequency droop gain.
The conventional 3rd-order model, \icl{differs from the detailed model (i.e. gives wrong stability results), for a large range of}
values of the droop gains.
\icl{It should also be noted} that the 5th-order and high-fidelity 3rd-order models show stability for only \icl{relatively a small range of values of the droop gains compared to those of the detailed average model.}

\paragraph{Eigenloci}	
\icl{These are shown in Fig. \ref{MV}
as the frequency droop coefficients $k_{pi}=k_{pk}$ vary in the interval  $[6\times10^{-5}, 5.3\times10^{-4}]$. The voltage droop coefficients are set to $k_{qi}= k_{qk}=1.5\times10^{-4}$} .
The eigenloci of the 5th-order and high-fidelity 3rd-order models are \icl{very close}, and show instability \icl{at lower values of the droop gains relative to those of the} 
detailed average model. This \icl{also} agrees with \icl{ Fig. \ref{regionMV}.} 
There is \icl{also} no significant change in the location of \icl{the} poles of the conventional 3rd-order model for \icl{large} droop gains, \icl{which} confirms its poor performance in \icl{Fig.~\ref{regionMV}.}

	\vspace{-3mm}
\begin{figure}[ht!]
	\begin{center}
		\includegraphics[width=0.4\textwidth]{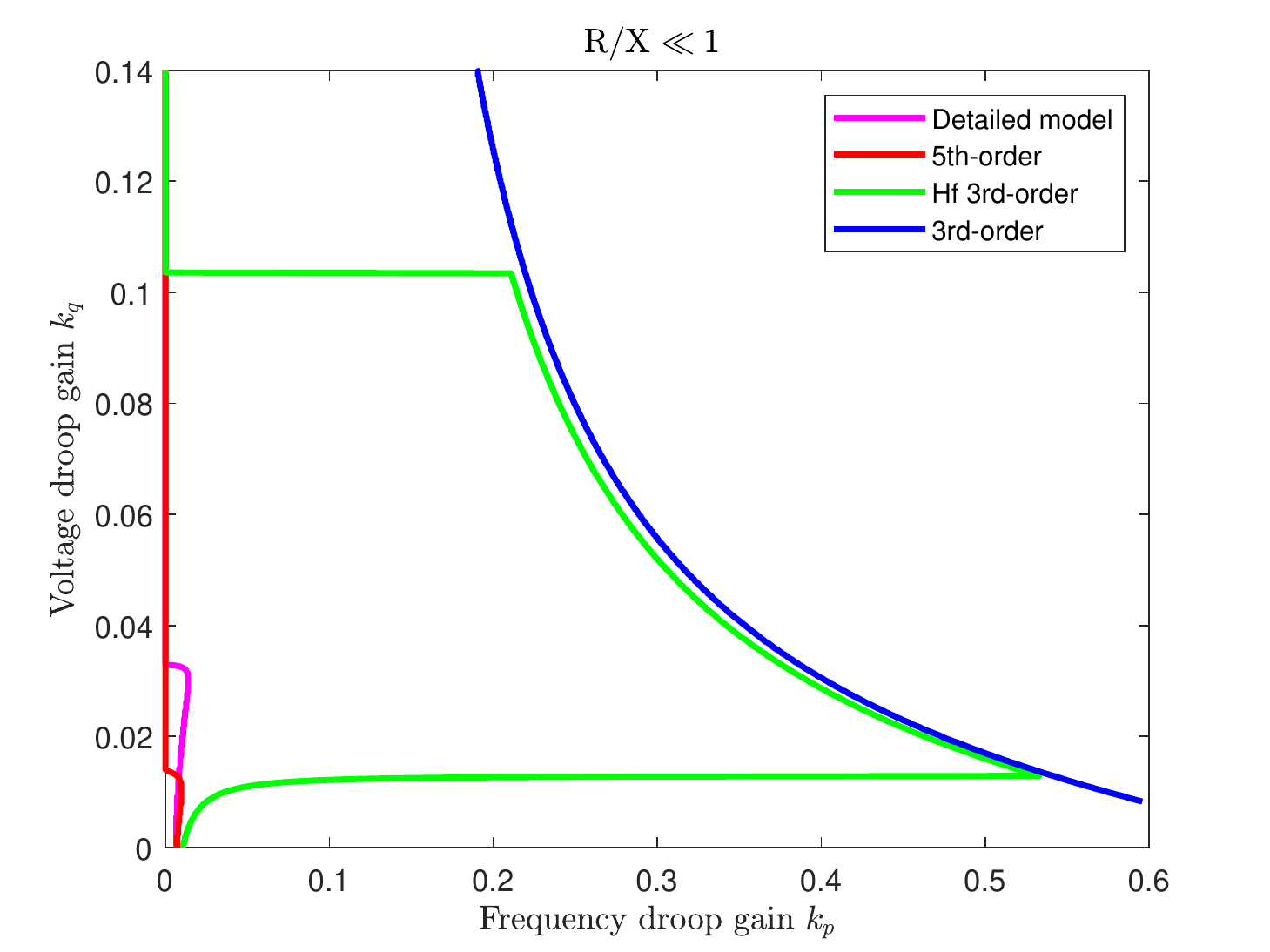}
		\vspace{-3mm}
		\caption{ Comparison of stability regions predicted by the four small-signal models for the case $R/X\ll1$. The region to the left of the corresponding curve are the values of the droop gains  for which the model is stable.}
		\vspace{-6mm}
		\label{regionHV}
	\end{center}
\end{figure}	

\begin{figure}[ht!]
	\begin{center}
		\includegraphics[width=0.4\textwidth]{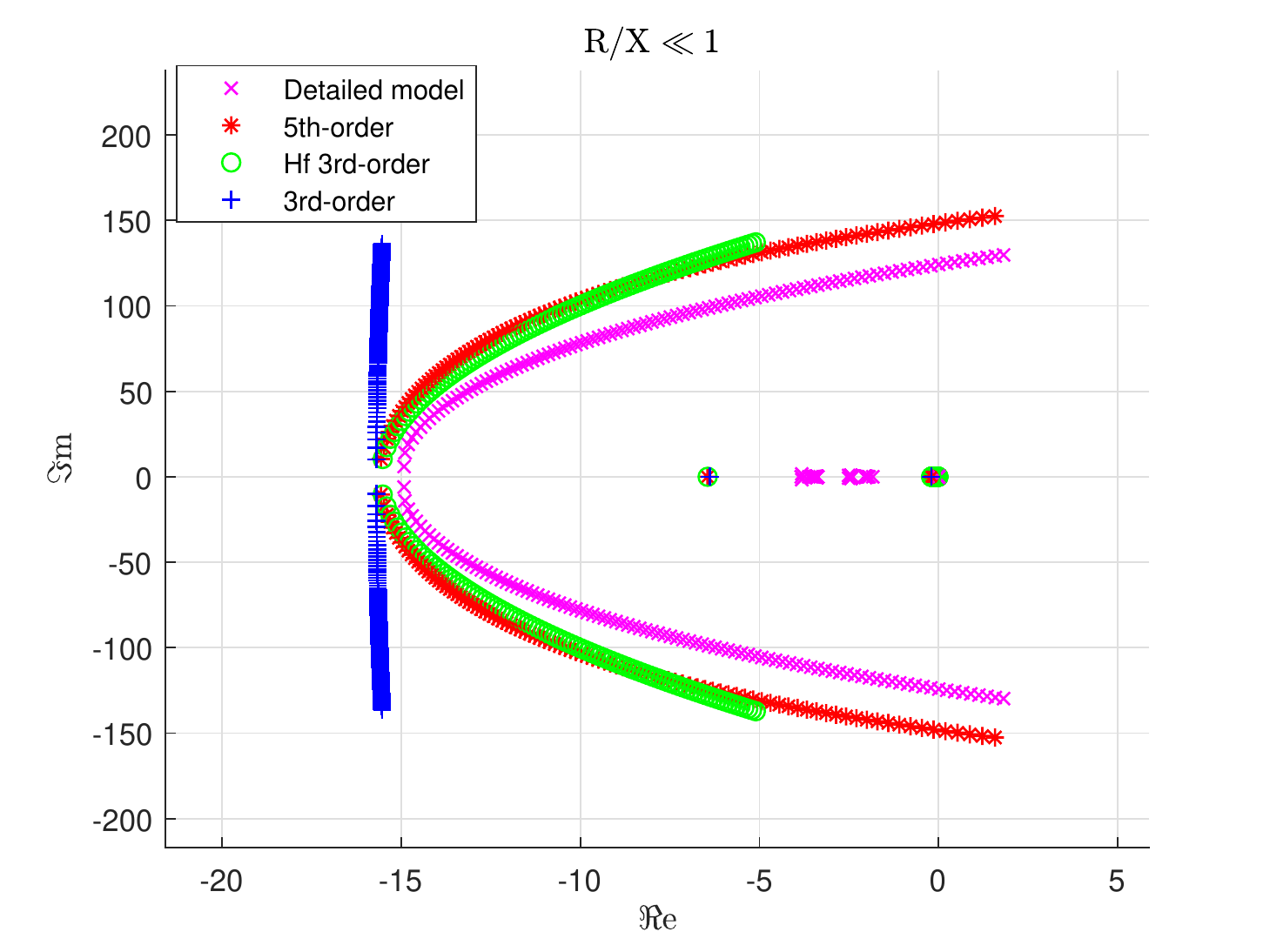}
		\vspace{-3mm}
		\caption{Eigenloci of the four small-signal models for the case $R/X\ll1$ when
			\icl{$k_{pi}=k_{pk}$ take values in the range $[6\times10^{-5}, 7\times10^{-3}]$ and
				$k_{qi}= k_{qk}=1.5\times10^{-4}$}.}
		\vspace{-5mm}
		\label{HV}
	\end{center}
\end{figure}

\subsubsection{Line ratio $R/X\ll1$}
\icl{To investigate this regime we considered}
line parameters \mbox{$R_{ik}=0.4~\Omega$,} \mbox{$L_{ik}=7~\text{mH}$}, \mbox{$R/X\approx0.18$} \cite{Guerrero2014, Guerrero2015, chandorkar1993}.

\paragraph{Stability Region}
\icl{These  are presented in Fig. \ref{regionHV}.}
	It is evident that the stability \icl{region} of the 5th-order model is reasonably within 
\icl{that} associated with the detailed average model, while those of the 3rd order models are erroneous for 
\icl{a large range of values of the} droop gains. This shows that the 3rd order models are highly unsuitable for stability assessment when $R/X\ll 1$.

	\paragraph{Eigenloci}
\icl{These are shown in Fig. \ref{HV}
as the frequency droop coefficients $k_{pi}=k_{pk}$ vary in the interval  $[6\times10^{-5}, 7\times10^{-3}]$. The voltage droop coefficients are set to $k_{qi}=k_{qk}=1.5\times10^{-4}$.}
The eigenloci of the 5th-order \icl{and the}  detailed average model shows instability
\icl{at lower values of the droop gains relative to those of the}
3rd order models.
	The location of \icl{the} poles of the conventional 3rd-order model does not change much for  large droop gains, \icl{and} hence incorrectly \icl{predicts stability in this regime as also shown}
in Fig. \ref{regionHV}. \icl{The pole movement of the high-fidelity 3rd-order model better reflects that of the detailed model, in comparison with the conventional 3rd-order model. There is, however, error in its stability predictions for large droop gains.}

\subsection{Nonlinear Model Dynamic Response}
\label{timedomain}
\begin{figure}[ht!]
	\begin{center}
		\includegraphics[width=0.5\textwidth]{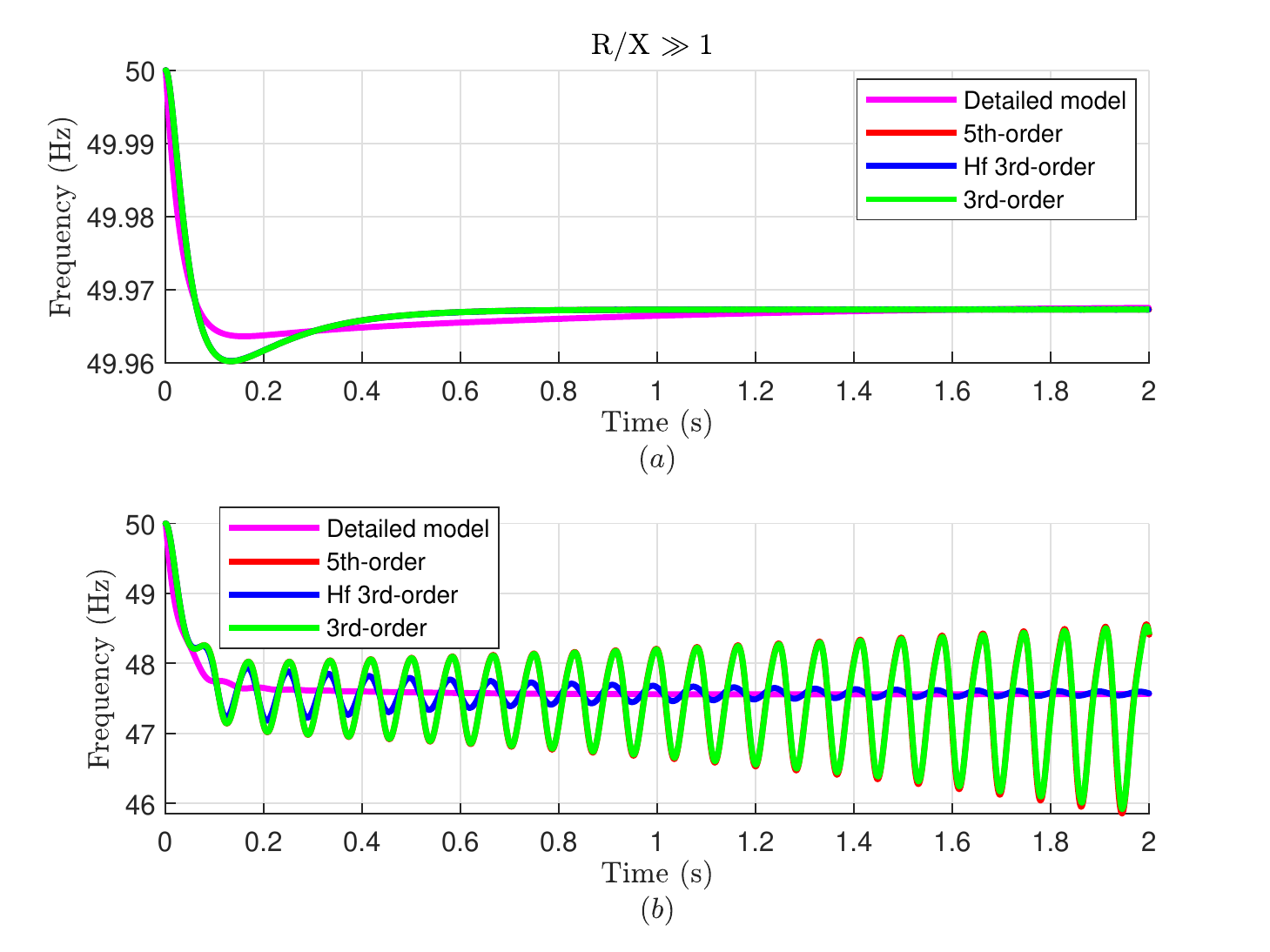}
		\vspace{-4mm}
	\caption{Dynamic responses of the four nonlinear models for the case $R/X\gg1$ when:  (a) \icl{$k_{pi}=k_{pk}=6\times10^{-5}$, $k_{qi}=k_{qk}=1.5\times10^{-4}$; (b) $k_{pi}=k_{pk}=4.4\times10^{-3}$, $k_{qi}=k_{qk}=1.5\times10^{-4}$}.}
		\rmspace
		\label{freqGG1}
	\end{center}
\end{figure}

\begin{figure}[ht!]
	\begin{center}
		\includegraphics[width=0.5\textwidth]{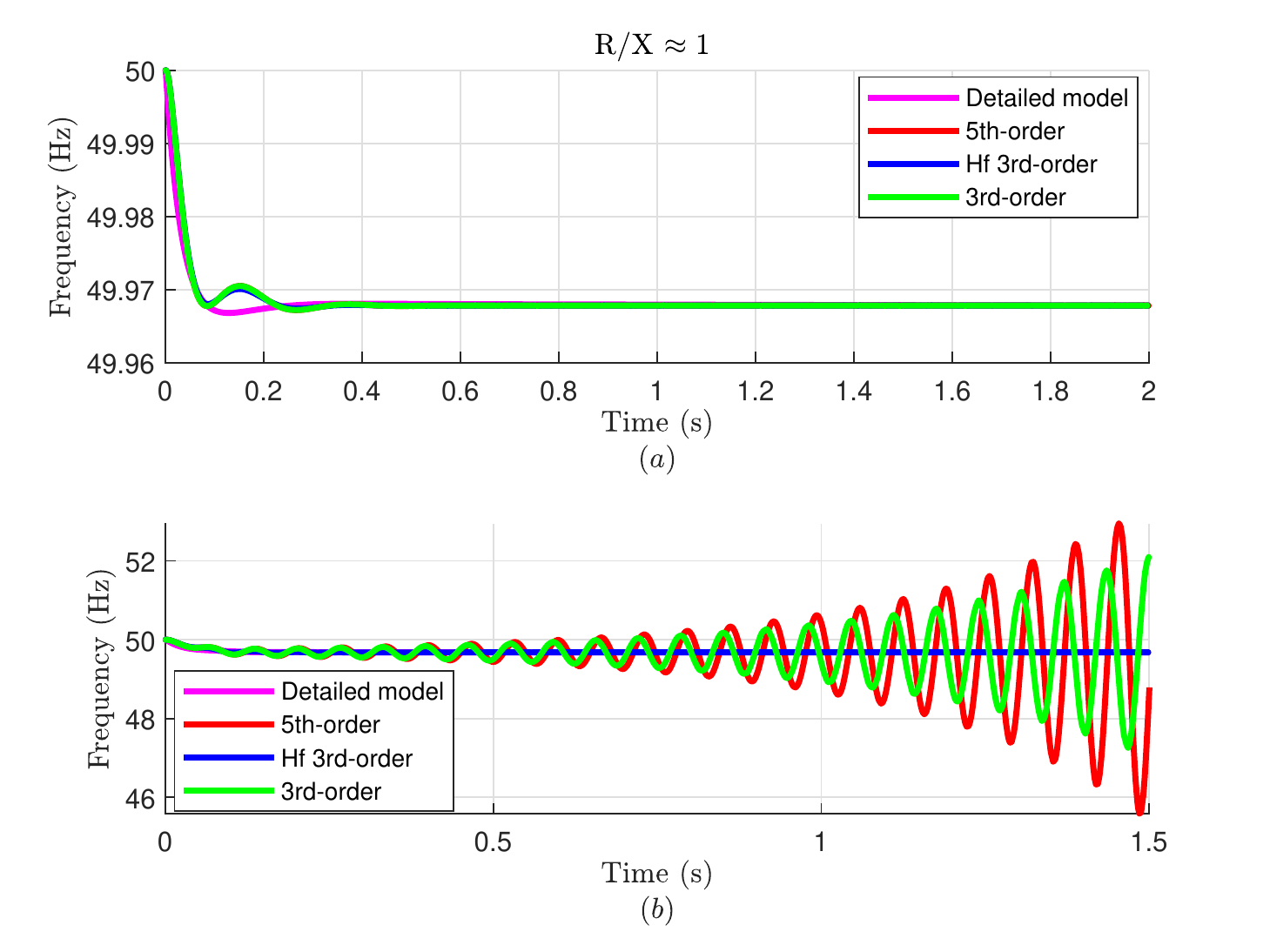}
		\vspace{-7mm}
		\caption{Dynamic responses of the four nonlinear models for the case $R/X\approx1$ when:  (a) \icl{$k_{pi}=k_{pk}=6\times10^{-5}$, $k_{qi}=k_{qk}=1.5\times10^{-4}$; (b) $k_{pi}=k_{pk}=5.3\times10^{-4}$, $k_{qi}=k_{qk}=1.5\times10^{-4}$}.}
		\vspace{-3mm}
		\label{freqEQ1}
	\end{center}
\end{figure}

\begin{figure}[ht!]
	\begin{center}
		\includegraphics[width=0.5\textwidth]{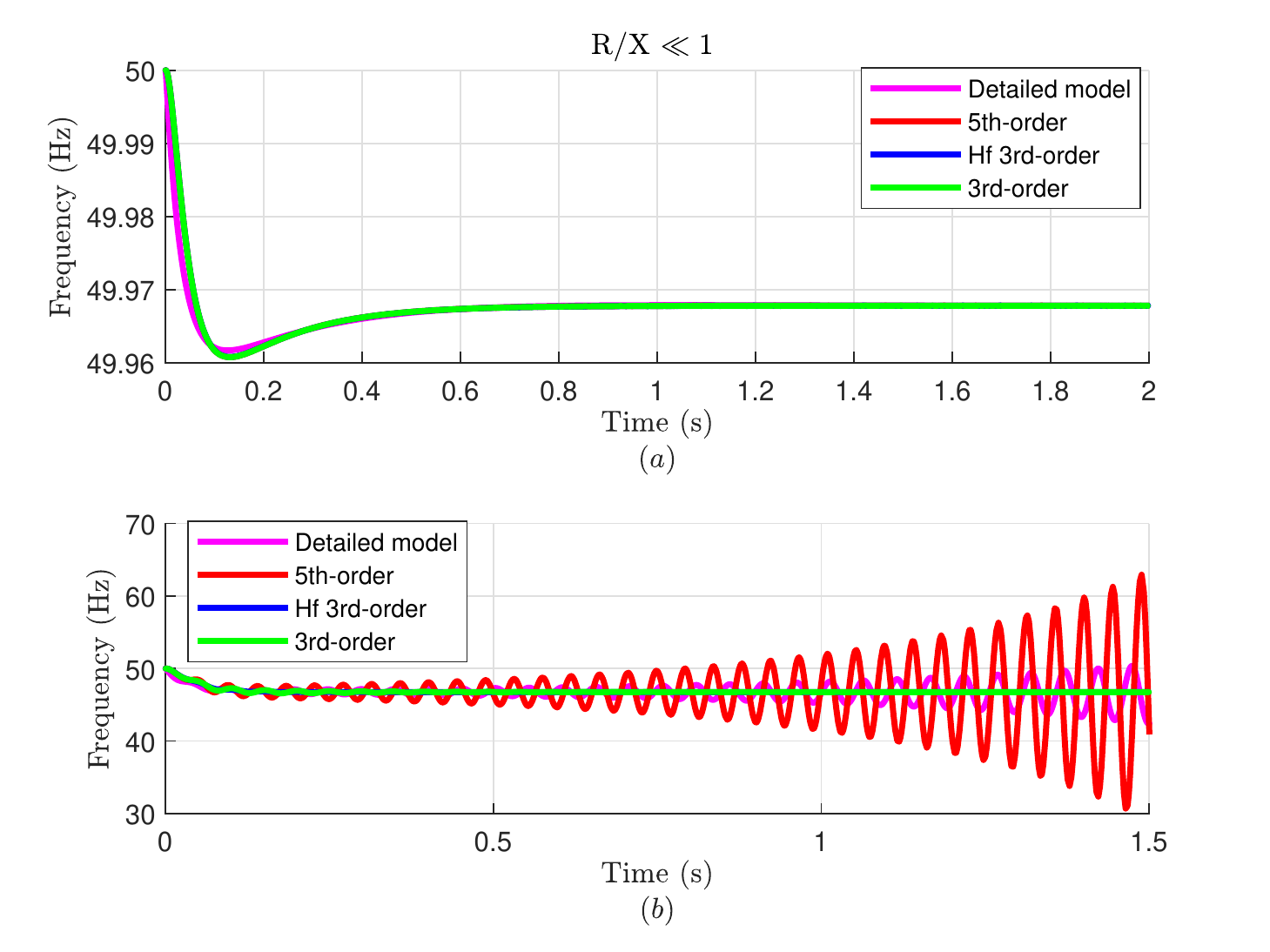}
		\vspace{-6mm}
		\caption{Dynamic responses of the four nonlinear models for the case $R/X\ll1$ when: (a) \icl{$k_{pi}=k_{pk}=6\times10^{-5}$, $k_{qi}=k_{qk}=1.5\times10^{-4}$; (b) $k_{pi}=k_{pk}=7\times10^{-3}$, $k_{qi}=k_{qk}=1.5\times10^{-4}$}.}
		\vspace{-6mm}
		\label{freqLL1}
	\end{center}
\end{figure}

In this section we \icl{present dynamic responses} that \icl{demonstrate} the stability results shown in Figs. \ref{LV}, \ref{MV}, \ref{HV} for the respective $R/X$ cases.
The \icl{model used} in the simulations \icl{is a more detailed one and}
includes the on/off actuation of the electronic switches \icl{via PWM.}
Figs. \ref{freqGG1}a, \ref{freqEQ1}a, \ref{freqLL1}a show the responses when the lower value of the frequency droop gain  \icl{used in} 
Figs. \ref{LV}, \ref{MV}, \ref{HV} is chosen, which \icl{is \icl{$k_{pi}=k_{pk}=6\times10^{-5}$}.}
Likewise Figs. \ref{freqGG1}b, \ref{freqEQ1}b, \ref{freqLL1}b show the responses when  the upper values of the frequency droop gain \icl{used in} 
Figs. \ref{LV}, \ref{MV}, \ref{HV} \icl{are} chosen, which are \icl{$k_{pi}=k_{pk}=4.4\times10^{-3}$, $5.3\times10^{-4}$}, $7\times10^{-3}$ \icl{respectively.}

It is evident in Figs. \ref{freqGG1}a, \ref{freqEQ1}a, \ref{freqLL1}a that the reduced-order models produce \icl{stable} responses \icl{which are similar to those of} the detailed average 
\icl{model when the lower values of the frequency droop gains are used}. This is expected as small droop gains do not cause instability \icl{(at the expense of a slower response)}.
The responses presented in Figs. \ref{freqGG1}b, \ref{freqEQ1}b, \ref{freqLL1}b, \icl{which demonstrate also unstable behaviour,} agree with the stability results in Figs. \ref{LV}, \ref{MV}, \ref{HV} respectively.
This further validates the veracity of the stability regions shown in Figs. \ref{regioLV}, \ref{regionMV}, \ref{regionHV}.

\section{Discussion}
\label{disc}
In this section we \icl{discuss the findings of the analysis in section  \ref{accuracy} and describe the relative merits of the models considered}.

The small-signal stability analysis showed that  the conventional 3rd-order model \icl{gives inaccurate stability results in all three regimes of the line $R/X$ ratio, and hence it cannot be recommended for stability analysis.}
The high-fidelity 3rd-order model performance is acceptable only when $R/X\approx1$, and
\icl{appears to be conservative, in the sense that it has a smaller stability region}
compared to that of the detailed average model.
The 5th-order model generally performs better and admits small error in its stability results  in all the cases except when $R/X\gg1$.


\yo{The \icl{discrepancy} observed in the stability results given by the reduced-order models \icl{relative} to \icl{those} of the detailed average model shows that neglecting the \icl{inverters'} current and voltage \icl{controller} dynamics affects their \icl{stability properties}. The additional \icl{omission} of the line dynamics in the conventional 3rd-order \icl{model} further explains its overall poor performance in all the \icl{cases,} compared to \icl{that} of the  other two reduced-order models. }
This shows \icl{that line dynamics can play a vital role in microgrid stability analysis} and  should not \icl{in general} be \icl{omitted.} 

\yo{The stablity \icl{of the} 
the detailed average model for a smaller range of the droop gains when $R/X\gg1$ suggests  that the inverter's voltage and current controllers, and \icl{their interaction with the line dynamics,} play a key role in the stability of very resistive microgrids. This \icl{is consistent with the insight} that when the $R/X$ value is large, the line dynamics become \icl{fast. Therefore} \icl{their timescale overlaps} that of the \icl{inverters'} voltage and current controllers, \icl{which  causes} their strong \icl{coupling.}}
Since the reduced models have poor performance when $R/X\gg1$, we recommend that the detailed average model is used in this case.


We also presented \icl{in section \ref{accuracy}} dynamic responses \icl{of the models considered and a comparison was also made with the responses of a more detailed switching model. The dynamic responses were consistent with the stability analysis and also validated the accuracy of the detailed average model.}

Based on our findings, Table \ref{table2} summarises the \icl{accuracy} 
of the three reduced-order models \icl{for predicting stability in} the three $R/X$ regimes considered.
We categorise this based on a descending scale of "Good", "Acceptable", and "Unacceptable".
 "Good" means that the \icl{stable region of the model is} within that of the detailed average model, even though some \icl{conservatism} may be present; "Acceptable" indicates that the associated stable region  is \icl{generally within} that of the detailed model, but \icl{there is a} small error for some droop gains;
	"Unacceptable" implies that \icl{erroneous stability results are} given for \icl{a large range of values of the} droop gains \icl{(i.e., the corresponding stable region is larger than that of the detailed average model)}.

	\begin{table}[ht!] 
		\caption{Stability Predictions of the Three Reduced-Order models
		}
		\vspace{-3mm}
		\small
		\label{table2}
		\begin{center}
			\begin{tabular}{cccc}
				\hline
				&$R/X\gg1$  &$R/X\approx1$ &$R/X\ll1$\\
				\hline
				EM 5th-order &Unacceptable  &Acceptable &Acceptable\\
				Hf 3rd-order &Unacceptable &Good &Unacceptable\\
				3rd-order  &Unacceptable  &Unacceptable &Unacceptable\\
				\hline
			\end{tabular}
		\end{center}
	\end{table}	

\subsection*{\textbf{Recommendations}}
It would be expected that the reduced-order models may admit some inaccuracies compared to the detailed average model in their \icl{stability properties}. \jdw{However, \icl{in order to avoid excessive discrepancies we present some brief recommendations on selecting an appropriate reduced order model, based on the results of our \icl{study}.}}

\begin{itemize}
	 \item The  high-fidelity 3rd-order and 5th-order models can be used in the case of \icl{line} $R/X$ ratio close to unity.
	 \item For inductive microgrids ($R/X\ll1$) we recommend that only the 5th-order model can be used.
	 \item We caution the use \icl{of reduced-order models that aim to simplify 
the inverter and line dynamics}
in the case of very resistive microgrids ($R/X\gg1$). 
\icl{The detailed average model is still though an appropriate model in this regime.}

	\end{itemize}

\section{Conclusion}
\label{concl}
\icl{Reduced order dynamic models for microgrids have been extensively used in the literature so as 
to facilitate system-wide stability assessment and  analytical studies.}
The 
\icl{simplifications} inevitably present in the reduced-order models raise the concern over the correctness of their stability \icl{properties} compared to \icl{those of more detailed average models}. Their performance also \icl{differs in different regimes} of the line $R/X$ ratios \icl{as the latter is associated with the timescale of the line dynamics and the extent to which these interfere with the dynamics of the inverters}.

We have therefore conducted a comprehensive comparative study of the accuracy of the \icl{following commonly used reduced order models:
the electromagnetic 5th-order, the conventional 3rd-order, and a high-fidelity 3rd-order model. A comparison has been made of the stability predictions of these models to those of a detailed average model for \icl{various regimes} of the line $R/X$ ratios. Our study has demonstrated that the 3rd-order model, where line dynamics are omited, can provide inaccurate stability results in all three regimes of the line $R/X$ ratios. The accuracy of the high-fidelity 3rd-order model
degrades as the $R/X$ ratio becomes either very large or small, and the electromagnetic 5th-order becomes inaccurate when the $R/X$ ratio is very large.}



\icl{Therefore the 5th-order and the high-fidelity 3rd-order models appear to be appropriate 
for microgrids with $R/X$ close to unity, and the former appears to be also quite accurate 
when the microgrid is inductive.}
We strongly advise caution when using reduced-order models in highly resistive microgrids. 

		

	
\appendices


\section{Detailed Average Model}  \label{ssdt}

In order to present the  small-signal linearized model, we introduce the error states: $\tilde\delta_i=\delta_i-\delta^s_i$, $\tilde\omega_i=\omega_i-\omega_0$, $\tilde V_i=V_i-V^s_i$,  $\tilde \phi_{dqi}=\phi_{dqi}-\phi^s_{dqi}$, $\tilde \gamma_{dqi}=\gamma_{dqi}-\gamma^s_{dqi}$, $\tilde v_{odqi}= v_{odqi}-v^s_{odqi}$,
$\tilde I_{DQ,ik}=I_{DQ,ik}-I^s_{DQ,ik}$,
$\tilde I_{DQ,\ell i}=I_{DQ,\ell i}-I^s_{DQ,\ell i}$.\newline

The entries of the Jacobian matrix for the detailed average model are defined as follows: 
{\small \mbox{Let  $r=i,k$,
		$\Lambda_p=\text{diag}\left(\frac{1}{k_{pi}}, \frac{1}{k_{pk}}\right)$,
		$\Lambda_q=\text{diag}\left(\frac{1}{k_{qi}}, \frac{1}{k_{qk}}\right)$},
	\mbox{$\sigma=\text{diag}(\mathbf{e},\mathbf{e})$, 
		$\sigma^v_{kp}=\text{diag}(K_{PVi}\mathbf{e},K_{PVk}\mathbf{e})$},
	\mbox{$\sigma^{cv}_{kp}=\text{diag}(K_{PVi}K_{PCi}\mathbf{e},K_{PVk}K_{PCk}\mathbf{e})$},
	\mbox{$K_{PV}=\text{diag}(K_{PVi}\mathbf{I}_2,K_{PVk}\mathbf{I}_2)$},
	\mbox{$K_{IV}=\text{diag}(K_{IVi}\mathbf{I}_2,K_{IVk}\mathbf{I}_2)$},
	\mbox{$K_{IC}=\text{diag}(K_{ICi}\mathbf{I}_2,K_{ICk}\mathbf{I}_2)$},
	$\sigma^v=\text{diag}(\sigma^v_i,\sigma^v_k), \sigma^v_r=K_{PCr}K_{PVr}\mathbf{I}_2+\mathbf{I}_2$, 	
	$\sigma_{pc}=\text{diag}(K_{PCi}K_{IVi}\mathbf{I}_2,K_{PCk}K_{IVk}\mathbf{I}_2)$,
	$\mathcal{Z}_{f}=\text{diag}(Z_{fi},Z_{fk})
	\in\mathbb{R}^{4\times4}$, 
	$Z_{fr}=-R_{fr}\mathbf{I}_2 + 2\omega_0 JL_{fr}-K_{PCr}\mathbf{I}_2$, 
	$Z_{fv}=\text{diag}(\omega_0C_{fi}J, \omega_0C_{fk}J)$,
	$\mathcal{Z}_{\ell}=\text{diag}(Z_i,Z_k)\in\mathbb{R}^{4\times4}$, $Z_r=-R_{\ell r}\mathbf{I}_2 + \omega_0 JL_{\ell r}, 
	\mathcal{Z}_{ik}=-R_{ik}\mathbf{I}_2 + \omega_0 JL_{ik}$,
	$N=\begin{bmatrix}	-\mathcal{T}(\delta^s_i) ~\mathcal{T}(\delta^s_k)
	\end{bmatrix}^{\top}$,
	\mbox{$N^{\delta}=\text{diag}(-\frac{\partial{\mathcal{T}^{-1}(\delta^s_i)}}{\partial\delta_i} I^s_{DQ,ik} , \frac{\partial{\mathcal{T}^{-1}(\delta^s_k)}}{\partial\delta_k} I^s_{DQ,ik} )$},
	\mbox{$N_{\ell}=\text{diag}(\mathcal{T}(\delta^s_i)^{-1},\mathcal{T}(\delta^s_k)^{-1})$},
	$C_f=\text{diag}(C_{fi}\mathbf{I}_2,C_{fk}\mathbf{I}_2)$,
	\mbox{$N^{\delta}_{\ell}=\text{diag}(\frac{\partial{\mathcal{T}^{-1}(\delta^s_i)}}{\partial\delta_i} I^s_{DQ\ell i} , \frac{\partial{\mathcal{T}^{-1}(\delta^s_k)}}{\partial\delta_k} I^s_{DQ\ell k} )$},
	\mbox{$\mathcal{N}^v=\text{diag}(\mathcal{T}(\delta^s_i) , \mathcal{T}(\delta^s_k))$},\\
	\mbox{$\mathcal{N}^{\delta}=\text{diag}(\frac{\partial{\mathcal{T}(\delta^s_i)}}{\partial\delta_i}v^s_{odqi}, \frac{\partial{\mathcal{T}(\delta^s_k)}}{\partial\delta_k}v^s_{odqk})$},
	\mbox{$\mathcal{D}^I_{\omega}=2\times\text{diag}(L_{fi}Ji^s_{dqi}, L_{fk}Ji^s_{dqk})$},
	\mbox{$L_{\ell}=\text{diag}(L_{\ell i}\mathbf{I}_2,L_{\ell k}\mathbf{I}_2)$},
	\mbox{$\mathcal{D}^v_{\omega}=\text{diag}(C_{fi}Jv^s_{odqi}, C_{fk}Jv^s_{odqk})$},
	\mbox{$\mathcal{\bar N}^{\delta}=\begin{bmatrix}\frac{\partial{\mathcal{T}(\delta^s_i)}}{\partial\delta_i}v^s_{odqi} ~- \frac{\partial{\mathcal{T}(\delta^s_k)}}{\partial\delta_k}v^s_{odqk} \end{bmatrix}$},
	\mbox{$\mathcal{N}=\begin{bmatrix}\mathcal{T}(\delta^s_i) ~- \mathcal{T}(\delta^s_k) \end{bmatrix}$},
	\mbox{$\frac{\partial P}{\partial v_{odq}}=\text{diag}(\frac{\partial P_i}{\partial v_{odqi}}\Bigr|_{x^s}, \frac{\partial P_k}{\partial v_{odqk}}\Bigr|_{x^s}), \bar N=N^{\delta}-N^{\delta}_{\ell}$}, 
	\mbox{$\frac{\partial Q}{\partial  v_{odq}}=\text{diag}(
		\frac{\partial Q_i}{\partial  v_{odqi}}\Bigr|_{x^s}, \frac{\partial Q_k}{\partial  v_{odqk}}\Bigr|_{x^s})$},
	$L_f=\text{diag}(L_{fi}\mathbf{I}_2,L_{fk}\mathbf{I}_2)$,
	\mbox{$\mathcal{D}^{q}_{\ell}=\text{diag}(D_i,D_k)\in\mathbb{R}^{2\times4},
		D_r=\begin{bmatrix}
		\frac{\partial Q_r}{\partial I_{D,\ell r}} \Bigr|_{x^s} & \frac{\partial Q_r}{\partial I_{Q,\ell r}} \Bigr|_{x^s}\end{bmatrix}$},
	\mbox{$\mathcal{D}^{p}_{\ell}=\text{diag}(\bar D_i,\bar D_k)\in\mathbb{R}^{2\times4}$, $\bar D_r=\begin{bmatrix}
		\frac{\partial P_r}{\partial I_{D,\ell r}} \Bigr|_{x^s} & \frac{\partial P_r}{\partial I_{Q,\ell r}} \Bigr|_{x^s} \end{bmatrix}$},
	\mbox{$\frac{\partial P}{\partial\delta}=\text{diag}\left(\frac{\partial P_i}{\partial\delta_i} \Bigr|_{x^s},  \frac{\partial P_k}{\partial\delta_k} \Bigr|_{x^s}\right)$},
	\mbox{$\frac{\partial Q}{\partial\delta}=\text{diag}\left(\frac{\partial Q_i}{\partial\delta_i} \Bigr|_{x^s},  \frac{\partial Q_k}{\partial\delta_k} \Bigr|_{x^s}\right)$},
	\mbox{$\frac{\partial P}{\partial I_{DQ,ik}}=\begin{bmatrix}
		\frac{\partial P_i}{\partial I_{D,ik}} \Bigr|_{x^s}& \frac{\partial P_i}{\partial I_{Q,ik}}\Bigr|_{x^s}\\\frac{\partial P_k}{\partial I_{D,ik}} \Bigr|_{x^s}& \frac{\partial P_k}{\partial I_{Q,ik}}\Bigr|_{x^s}
		\end{bmatrix}$},
	\mbox{$\frac{\partial Q}{\partial I_{DQ,ik}}=\begin{bmatrix}
		\frac{\partial Q_i}{\partial I_{D,ik}} \Bigr|_{x^s}& \frac{\partial Q_i}{\partial I_{Q,ik}}\Bigr|_{x^s}\\\frac{\partial Q_k}{\partial I_{D,ik}} \Bigr|_{x^s}& \frac{\partial Q_k}{\partial I_{Q,ik}}\Bigr|_{x^s}
		\end{bmatrix}$},\\
	$\Gamma=\text{diag}(\mathbf{I}_{2}, \tau\Lambda_p, \tau\Lambda_q, \mathbf{I}_{4}, \mathbf{I}_{4}, L_{f}, C_{f}, L_{ik}\mathbf{I}_{2},  L_{\ell })\in\mathbb{R}^{28\times28}$.
}

The small-signal state space model of the detailed average model is described by \eqref{statedt}:
\begin{equation}
\label{statedt}
\dot{\tilde{x}}=\Gamma^{-1}\begin{bmatrix}
A_{11} & A_{12} \\ A_{21} & A_{22}
\end{bmatrix}\tilde x
\end{equation}
where 
\begin{align*} 
\tilde{x}=[\tilde\delta_i, \tilde\delta_k,\tilde\omega_i,\tilde\omega_k,
\tilde \phi^{\top}_{dqi}, \tilde \phi^{\top}_{dqk}, \tilde\gamma^{\top}_{dqi}, \tilde\gamma^{\top}_{dqk},
\tilde i^{\top}_{dqi}, \tilde i^{\top}_{dqk},\\
\tilde v^{\top}_{odqi}, \tilde v^{\top}_{odqk},
\tilde I^{\top}_{DQ,ik}, \tilde I^{\top}_{DQ,\ell i}, \tilde I^{\top}_{DQ,\ell k}]^{\top}
\end{align*}
{\scriptsize
	$ A_{11}=\begin{bmatrix}
	0  &\mathbf{I}_2 &0 &0 &0 &0 \\
	-\frac{\partial P}{\partial\delta} &-\Lambda_p &0 &0 &0 &0 \\
	-\frac{\partial Q}{\partial\delta} &0 &-\Lambda_q &0 &0 &0  \\
	0 &0 &\sigma &0 &0 &0 \\
	0 &0 &\sigma^{v}_{kp} &K_{IV} &0 &-\mathbf{I}_{4} \\
	0 &\mathcal{D}^I_{\omega} &\sigma^{cv}_{kp} &\sigma_{pc} &K_{IC} &Z_f
	\end{bmatrix}$\\
	$ A_{12}=\begin{bmatrix}
	0 &0 &0 \\  -\frac{\partial P}{\partial v_{odq}} &-\frac{\partial P}{\partial I_{DQ,ik}} &-\mathcal{D}^{p}_{\ell}\\ -\frac{\partial Q}{\partial v_{odq}} &-\frac{\partial Q}{\partial I_{DQ,ik}} &-\mathcal{D}^{q}_{\ell} \\ -\mathbf{I}_{4} &0 &0 \\ -K_{PV} &0 &0 \\ -\sigma^v &0 &0 \\ 
	\end{bmatrix}$\\
	$A_{21}=\begin{bmatrix}
	\bar N &\mathcal{D}^v_{\omega} &0 &0 &0 &\mathbf{I}_4 \\
	\mathcal{\bar N}^{\delta} &0 &0 &0 &0 &0 \\
	\mathcal{N}^{\delta} &0 &0 &0 &0 &0 
	\end{bmatrix}$\\
	$A_{22}=\begin{bmatrix}
	Z_{fv}  &N &-N_{\ell} \\ \mathcal{N} &\mathcal{Z}_{ik} &0 \\ \mathcal{N}^{v} &0 &\mathcal{Z}_{\ell}
	\end{bmatrix}$\\
}
$\tilde{x}\in\mathbb{R}^{28}, A\in\mathbb{R}^{28\times28}$.

\section{EM 5th-Order Small-Signal Model}	\label{ss5}
The entries of the Jacobian matrix for the EM 5th-order model are defined as follows: \\{\small
	$\mathcal{M}^{\delta}_{\ell}=\text{diag}(\frac{\partial{\mathcal{T}(\delta^s_i)}}{\partial\delta_i}\mathbf{e}V^s_i, \frac{\partial{\mathcal{T}(\delta^s_k)}}{\partial\delta_k}\mathbf{e}V^s_k)$,\\ 
	$\mathcal{M}^{v}=[\mathcal{T}(\delta_i)\mathbf{e}~ -\mathcal{T}(\delta_k)\mathbf{e}]$,\\
	$\mathcal{M}^{v}_{\ell}=\text{diag}(\mathcal{T}(\delta_i)\mathbf{e}, \mathcal{T}(\delta_k)\mathbf{e})$, \\
	$\mathcal{M}^{\delta}=[\frac{\partial{\mathcal{T}(\delta^s_i)}}{\partial\delta_i}\mathbf{e}V^s_i ~-\frac{\partial{\mathcal{T}(\delta^s_k)}}{\partial\delta_k}\mathbf{e}V^s_k  ]$,
	$\frac{\partial P}{\partial V}=\text{diag}(
	\frac{\partial P_i}{\partial V_i}\Bigr|_{x^s}, \frac{\partial P_k}{\partial V_k}\Bigr|_{x^s})$, 
	$\frac{\partial Q}{\partial V}=\text{diag}(
	\frac{\partial Q_i}{\partial V_i}\Bigr|_{x^s}, \frac{\partial Q_k}{\partial V_k}\Bigr|_{x^s})$,\\
	$\Gamma_5=\text{diag}(\mathbf{I}_{2}, \tau\Lambda_p, \tau\Lambda_q, L_{ik}\mathbf{I}_{2}, L_{\ell })\in\mathbb{R}^{12\times12}$.
}

The 5th-order small-signal state space model is described by \eqref{state5}:
\begin{equation}
\label{state5}
\dot{\tilde{x}}_{5}=\Gamma_5^{-1}A_{5}\tilde x_{5}
\end{equation}
where 
$\tilde{x}_{5}=[\tilde\delta_i, \tilde\delta_k,\tilde\omega_i,\tilde\omega_k,\tilde V_i, \tilde V_k, \tilde I^{\top}_{DQ,ik}, \tilde I^{\top}_{DQ,\ell i}, \tilde I^{\top}_{DQ,\ell k}]^{\top},$
\begin{equation*}\scriptsize
A_{5}=\begin{bmatrix}
\mymathbb{0}_{2\times2} &\mathbf{I}_{2} &\mymathbb{0}_{2\times2} &\mymathbb{0}_{2\times2} &\mymathbb{0}_{2\times4} \\
-\frac{\partial P}{\partial\delta} &-\Lambda_p &-\frac{\partial P}{\partial V} &-\frac{\partial P}{\partial I_{dq,ik}} &-\mathcal{D}^{p}_{\ell}\\
-\frac{\partial Q}{\partial\delta} &\mymathbb{0}_{2\times2} &-(\Lambda_q+\frac{\partial Q}{\partial V}) &-\frac{\partial Q}{\partial I_{dq,ik}} &-\mathcal{D}^{q}_{\ell}\\
\mathcal{M}^{\delta} &\mymathbb{0}_{2\times2} &\mathcal{M}^{v} &\mathcal{Z}_{ik} &\mymathbb{0}_{2\times4}\\
\mathcal{M}^{\delta}_{\ell} &\mymathbb{0}_{4\times2} &\mathcal{M}^{v}_{\ell} &\mymathbb{0}_{4\times2 } & \mathcal{Z}_{\ell}
\end{bmatrix},
\end{equation*}
$\tilde{x}_{5}\in\mathbb{R}^{12}, A_{5}\in\mathbb{R}^{12\times12}$. Other entries in $\Gamma_{5}, A_5$ are as defined in Appendix \ref{ssdt}

\section{Conventional 3rd-Order Small-Signal Model}\label{ss3}
The following are defined  to present the Jacobian matrix of the 3rd-order small-signal model:\\
{ \scriptsize $\mathcal{D}^p_{\delta}=\begin{bmatrix}
	\frac{\partial P_i}{\partial \delta_i}\Bigr|_{x^s}
	&\frac{\partial P_i}{\partial \delta_k}\Bigr|_{x^s}\\
	\frac{\partial P_k}{\partial \delta_i}\Bigr|_{x^s}
	&\frac{\partial P_k}{\partial \delta_k}\Bigr|_{x^s}
	\end{bmatrix}$,
	$\mathcal{D}^q_{\delta}=\begin{bmatrix}
	\frac{\partial Q_i}{\partial \delta_i}\Bigr|_{x^s}
	&\frac{\partial Q_i}{\partial \delta_k}\Bigr|_{x^s}\\
	\frac{\partial Q_k}{\partial \delta_i}\Bigr|_{x^s}
	&\frac{\partial Q_k}{\partial \delta_k}\Bigr|_{x^s}
	\end{bmatrix}$,\\
	$\mathcal{D}^p_{v}=\begin{bmatrix}
	\frac{\partial P_i}{\partial V_i}\Bigr|_{x^s}
	&\frac{\partial P_i}{\partial V_k}\Bigr|_{x^s}\\
	\frac{\partial P_k}{\partial V_i}\Bigr|_{x^s}
	&\frac{\partial P_k}{\partial V_k}\Bigr|_{x^s}
	\end{bmatrix}$, 
	$\mathcal{D}^q_{v}=\begin{bmatrix}
	\frac{\partial Q_i}{\partial V_i}\Bigr|_{x^s}
	&\frac{\partial Q_i}{\partial V_k}\Bigr|_{x^s}\\
	\frac{\partial Q_k}{\partial V_i}\Bigr|_{x^s}
	&\frac{\partial Q_k}{\partial V_k}\Bigr|_{x^s}
	\end{bmatrix}$,\\
	$\Gamma_3=\text{diag}(\mathbf{I}_{2}, \tau\Lambda_p, \tau\Lambda_q, L_{\ell })\in\mathbb{R}^{10\times10}$}.

The conventional 3th-order small-signal state space model is described by \eqref{state3}:
\begin{equation}
\label{state3}
\dot{\tilde{x}}_{3}=\Gamma_3^{-1}A_{3}\tilde x_{3}
\end{equation}
where 
$\tilde{x}_{3}=[\tilde\delta_i, \tilde\delta_k,\tilde\omega_i,\tilde\omega_k,\tilde V_i, \tilde V_k,  \tilde I^{\top}_{dq,\ell i}, \tilde I^{\top}_{dq,\ell k}]^{\top},$
\begin{equation*} \scriptsize
A_3=\begin{bmatrix}
\mymathbb{0}_{2\times2} &\mathbf{I}_{2}  &\mymathbb{0}_{2\times2} &\mymathbb{0}_{2\times4} \\
-\mathcal{D}^p_{\delta}  &-\Lambda_p  &-\mathcal{D}^p_{v}  &-\mathcal{D}^p_{\ell}\\
-\mathcal{D}^q_{\delta} &\mymathbb{0}_{2\times2} &-(\Lambda_q+\mathcal{D}^q_{v})  &-\mathcal{D}^q_{\ell}\\
\mathcal{M}^{\delta}_{\ell} &\mymathbb{0}_{4\times2} & \mathcal{M}^{v}_{\ell} & \mathcal{Z}_{\ell}
\end{bmatrix},
\end{equation*}
$\tilde{x}_{3}\in\mathbb{R}^{10}, A_{3}\in\mathbb{R}^{10\times10}$. Other entries in $\Gamma_{\text{3}}, A_3$ are as defined in Appendixes \ref{ssdt} and \ref{ss5} respectively.

\section{High-Fidelity 3rd-Order Small-Signal Model} \label{sshf}
We define the following to present the Jacobian matrix of the high-fidelity 3rd-order small-signal model:\\
{ \scriptsize $\mathcal{D}^p_{\omega}=\begin{bmatrix}
	\frac{\partial P_i}{\partial \dot\delta_i}\Bigr|_{x^s}
	&\frac{\partial P_i}{\partial \dot\delta_k}\Bigr|_{x^s}\\
	\frac{\partial P_k}{\partial \dot\delta_i}\Bigr|_{x^s}
	&\frac{\partial P_k}{\partial \dot\delta_k}\Bigr|_{x^s}
	\end{bmatrix}$,
	$\mathcal{D}^q_{\omega}=\begin{bmatrix}
	\frac{\partial Q_i}{\partial \dot\delta_i}\Bigr|_{x^s}
	&\frac{\partial Q_i}{\partial \dot\delta_k}\Bigr|_{x^s}\\
	\frac{\partial Q_k}{\partial \dot\delta_i}\Bigr|_{x^s}
	&\frac{\partial Q_k}{\partial \dot\delta_k}\Bigr|_{x^s}
	\end{bmatrix}$,\\
	$\mathcal{\bar D}^p_{v}=\begin{bmatrix}
	\frac{\partial P_i}{\partial \dot V_i}\Bigr|_{x^s}
	&\frac{\partial P_i}{\partial \dot V_k}\Bigr|_{x^s}\\
	\frac{\partial P_k}{\partial\dot V_i}\Bigr|_{x^s}
	&\frac{\partial P_k}{\partial\dot V_k}\Bigr|_{x^s} 	\end{bmatrix}$,
	$\mathcal{\bar D}^q_{v}=\begin{bmatrix}
	\frac{\partial Q_i}{\partial \dot V_i}\Bigr|_{x^s}
	&\frac{\partial Q_i}{\partial \dot V_k}\Bigr|_{x^s}\\
	\frac{\partial Q_k}{\partial\dot V_i}\Bigr|_{x^s}
	&\frac{\partial Q_k}{\partial\dot V_k}\Bigr|_{x^s}
	\end{bmatrix}$,\\
	\\
	$\Gamma_{\text{hf}}=\begin{bmatrix}
	\mathbf{I}_{2} &\mymathbb{0}_{2\times2} &\mymathbb{0}_{2\times2} &\mymathbb{0}_{2\times4}\\
	\mathcal{D}^p_{\omega} &\tau\Lambda_p  &\mathcal{\bar D}^p_{v} &\mymathbb{0}_{2\times4}\\
	\mathcal{D}^q_{\omega} &\mymathbb{0}_{2\times2} &\tau\Lambda_q+\mathcal{\bar D}^q_{v}  &\mymathbb{0}_{2\times4}\\
	\mymathbb{0}_{4\times2} &	\mymathbb{0}_{4\times2} &	\mymathbb{0}_{4\times2} &L_{\ell}
	\end{bmatrix}$.
}

The high-fidelity 3th-order small-signal state space model is described by \eqref{statehf}.
\begin{equation}
\label{statehf}
\dot{\tilde{x}}_{\text{hf}}=\Gamma_{\text{hf}}^{-1}A_{\text{hf}}\tilde x_{\text{hf}}
\end{equation}
where 
$\tilde{x}_{\text{hf}}=[\tilde\delta_i, \tilde\delta_k,\tilde\omega_i,\tilde\omega_k,\tilde V_i, \tilde V_k,  \tilde I^{\top}_{DQ,\ell i}, \tilde I^{\top}_{DQ,\ell k}]^{\top},$
\begin{equation*} \scriptsize
A_{\text{hf}}=\begin{bmatrix}
\mymathbb{0}_{2\times2} &\mathbf{I}_{2}  &\mymathbb{0}_{2\times2} &\mymathbb{0}_{2\times4} \\
-\mathcal{D}^p_{\delta}  &-\Lambda_p  &-\mathcal{D}^p_{v}  &-\mathcal{D}^p_{\ell}\\
-\mathcal{D}^q_{\delta} &\mymathbb{0}_{2\times2} &-(\Lambda_q+\mathcal{D}^q_{v})  &-\mathcal{D}^q_{\ell}\\
\mathcal{M}^{\delta}_{\ell} &\mymathbb{0}_{4\times2} & \mathcal{M}^{v}_{\ell} & \mathcal{Z}_{\ell}
\end{bmatrix},
\end{equation*}
$\tilde{x}_{\text{hf}}\in\mathbb{R}^{10}, A_{\text{hf}}\in\mathbb{R}^{10\times10}$. Other entries in $\Gamma_{\text{hf}}, A_{\text{hf}}, $ are as defined in Appendixes \ref{ssdt}, \ref{ss5} and \ref{ss3} respectively.

\end{document}